\documentstyle[12pt,bezier]{article}
\parskip=\medskipamount                 % space between paragraphs
\thicklines                     % thick straight lines for pictures
\newcommand{\bimn}[7]{\bibitem{#1}#2,
{\em #3},
{ #4}$\;${\bf
#5}$\;$(#6)$\;${#7}.}

\def\inbar{\vrule height1.5ex width.4pt depth0pt}
\def\IC{\relax\,\hbox{$\inbar\kern-.3em{\rm C}$}}
\def\IN{\relax{\rm I\kern-.18em N}}
\def\IQ{\relax\,\hbox{$\inbar\kern-.3em{\rm Q}$}}
\def\IR{\relax{\rm I\kern-.18em R}}
\def\ZZ{\relax{\sf Z\kern-.4em Z}}
\def\a{\alpha}    \def\e{\epsilon}

 \def\cB{{\cal B}} \def\cC{{\cal C}} 
   
 \def\cK{{\cal K}}  
 \def\cO{{\cal O}}  
\def\cR{{\cal R}} 
\newtheorem{theorem}{Theorem}[section]
\newtheorem{proposition}{Proposition}[section]
\newtheorem{corollary}{Corollary}[section]

\newtheorem{lemma}{Lemma}[section]

\catcode`\@=11

\newif\if@fewtab\@fewtabtrue

\catcode`\@=11

\newif\if@fewtab\@fewtabtrue

{\count255=\time\divide\count255 by 60
\xdef\hourmin{\number\count255}
\multiply\count255 by-60\advance\count255 by\time
\xdef\hourmin{\hourmin:\ifnum\count255<10 0\fi\the\count255}}
\def\ps@draft{\let\@mkboth\@gobbletwo
    \def\@oddhead{}
    \def\@oddfoot
      {\hbox to 7 cm{\footnotesize {\em Draft of \jobname:} \draftdate
       \hfil}\hskip -7cm\hfil\rm\thepage \hfil}
    \def\@evenhead{}\let\@evenfoot\@oddfoot}

\def\ceqno{\global\@fewtabfalse
    \ifcase\@eqcnt \def\@tempa{& & &}\or \def\@tempa{& &}
      \or \def\@tempa{&}
      \or\def\@tempa{}\fi\@tempa
{\rm(\theequation)}}

\def\aeqno#1{\global\@fewtabfalse
    \ifcase\@eqcnt \def\@tempa{& & &}\or \def\@tempa{& &}
      \or \def\@tempa{&}
      \or\def\@tempa{}\fi\@tempa
{\rm(\theequation,#1)}}

\def\label#1{\ifnum\draftcontrol=1
 \global\def\draftnote{$\scriptstyle #1$}\fi
 \@bsphack\if@filesw {\let\thepage\relax
   \def\protect{\noexpand\noexpand\noexpand}%
\xdef\@gtempa{\write\@auxout{\string
      \newlabel{#1}{{\@currentlabel}{\thepage}}}}}\@gtempa
   \if@nobreak \ifvmode\nobreak\fi\fi\fi
  \@esphack}

\def\alabel#1#2{\label{#1}\global\@fewtabfalse
    \ifcase\@eqcnt \def\@tempa{& & &}\or \def\@tempa{& &}
      \or \def\@tempa{&}
      \or\def\@tempa{}\fi\@tempa
{\hbox to 3cm{\phantom{\rm(\theequation,#2)}
\draftnote \hfil}\hskip -3cm {\rm(\theequation,#2)}}}

\def\clabel#1{\label{#1}\global\@fewtabfalse
    \ifcase\@eqcnt \def\@tempa{& & &}\or \def\@tempa{& &}
      \or \def\@tempa{&}
      \or\def\@tempa{}\fi\@tempa
{\hbox to 3cm{\phantom{\rm(\theequation)}
\draftnote \hfil}\hskip -3cm{\rm(\theequation)}}}

\def\eqnarray{\def\draftnote{{}}\global\@fewtabtrue
\stepcounter{equation}\let\@currentlabel=\theequation
\global\@eqnswtrue
\global\@eqcnt\z@\tabskip\@centering\let\\=\@eqncr
$$\halign to \displaywidth\bgroup\@eqnsel\hskip\@centering\@eqcnt\z@
  $\displaystyle\tabskip\z@{##}$&\global\@eqcnt\@ne
  \hskip 1\arraycolsep \hfil$\displaystyle{##}$\hfil
  &\global\@eqcnt\tw@ \hskip 1\arraycolsep
$\displaystyle\tabskip\z@{##}$
\hfil  \tabskip\@centering&\global\@eqcnt\thr@@\llap{##}\tabskip\z@
\cr}

\def\endeqnarray{\@@eqncr\egroup
      \global\advance\c@equation\m@ne$$\global\@ignoretrue}

\def\@eqnnum{\hbox to 3cm{\phantom{\rm(\theequation)} \draftnote
                         \hfil}\hskip -3cm {\rm(\theequation)}}

\def\@@eqncr{\let\@tempa\relax
    \ifcase\@eqcnt \def\@tempa{& & &}\or \def\@tempa{& &}
      \or \def\@tempa{&}
      \or\def\@tempa{}
\fi\@tempa
\if@eqnsw
\if@fewtab\@eqnnum\fi
\stepcounter{equation}\fi\global
\@eqnswtrue\global\@eqcnt\z@\global\@fewtabtrue\cr}

\def\draftcite#1{\ifnum\draftcontrol=1#1\else{}\fi}

\def\@lbibitem[#1]#2{\item{}\hskip -3cm \hbox to 2cm
{\hfil$\scriptstyle\draftcite{#2}$}\hskip
1cm[\@biblabel{#1}]\if@filesw
     {\def\protect##1{\string ##1\space}\immediate
      \write\@auxout{\string\bibcite{#2}{#1}}}\fi\ignorespaces}

\def\@bibitem#1{\item\hskip -3cm \hbox to 2cm
{\hfil $\scriptstyle\draftcite{#1}$}\hskip 1cm
\if@filesw \immediate\write\@auxout
       {\string\bibcite{#1}{\the\value{\@listctr}}}\fi\ignorespaces}

\def\nsection#1{\section{#1}\setcounter{equation}{0}}

\def\nappendixe{\def\thesection{A}\section*{Appendix }
\def\theequation{{A.\arabic{equation}}}
\def\theproposition{{A.\arabic{proposition}}}
\setcounter{equation}{0}
\setcounter{proposition}{0}}

\def\draftdate{\number\month/\number\day/\number\year\ \ \ \hourmin }

\global\def\draftcontrol{0}
\catcode`\@=12

\def\theequation{{\thesection.\arabic{equation}}}

\newcommand{\qbezier}{\bezier{120}}
\def\DottedCircle{
\bezier{4}(0.966,-0.259)(1.04,0)(0.966,0.259)
\bezier{4}(0.966,0.259)(0.897,0.518)(0.707,0.707)
\bezier{4}(0.707,0.707)(0.518,0.897)(0.259,0.966)
\bezier{4}(0.259,0.966)(0,1.04)(-0.259,0.966)
\bezier{4}(-0.259,0.966)(-0.518,0.897)(-0.707,0.707)
\bezier{4}(-0.707,0.707)(-0.897,0.518)(-0.966,0.259)
\bezier{4}(-0.966,0.259)(-1.04,0)(-0.966,-0.259)
\bezier{4}(-0.966,-0.259)(-0.897,-0.518)(-0.707,-0.707)
\bezier{4}(-0.707,-0.707)(-0.518,-0.897)(-0.259,-0.966)
\bezier{4}(-0.259,-0.966)(0,-1.04)(0.259,-0.966)
\bezier{4}(0.259,-0.966)(0.518,-0.897)(0.707,-0.707)
\bezier{4}(0.707,-0.707)(0.897,-0.518)(0.966,-0.259)
}
\def\Endpoint[#1]{
\ifcase#1
\put(1,0){\circle*{0.15}}
\or\put(0.866,0.5){\circle*{0.15}}
\or\put(0.5,0.866){\circle*{0.15}}
\or\put(0,1){\circle*{0.15}}
\or\put(-0.5,0.866){\circle*{0.15}}
\or\put(-0.866,0.5){\circle*{0.15}}
\or\put(-1,0){\circle*{0.15}}
\or\put(-0.866,-0.5){\circle*{0.15}}
\or\put(-0.5,-0.866){\circle*{0.15}}
\or\put(0,-1){\circle*{0.15}}
\or\put(0.5,-0.866){\circle*{0.15}}
\or\put(0.866,-0.5){\circle*{0.15}}
\fi}
\def\Arc[#1]{
\thicklines			% this can be changed!
\ifcase#1
\bezier{25}(0.966,-0.259)(1.04,0)(0.966,0.259)
\or
\bezier{25}(0.966,0.259)(0.897,0.518)(0.707,0.707)
\or
\bezier{25}(0.707,0.707)(0.518,0.897)(0.259,0.966)
\or
\bezier{25}(0.259,0.966)(0,1.04)(-0.259,0.966)
\or
\bezier{25}(-0.259,0.966)(-0.518,0.897)(-0.707,0.707)
\or
\bezier{25}(-0.707,0.707)(-0.897,0.518)(-0.966,0.259)
\or
\bezier{25}(-0.966,0.259)(-1.04,0)(-0.966,-0.259)
\or
\bezier{25}(-0.966,-0.259)(-0.897,-0.518)(-0.707,-0.707)
\or
\bezier{25}(-0.707,-0.707)(-0.518,-0.897)(-0.259,-0.966)
\or
\bezier{25}(-0.259,-0.966)(0,-1.04)(0.259,-0.966)
\or
\bezier{25}(0.259,-0.966)(0.518,-0.897)(0.707,-0.707)
\or
\bezier{25}(0.707,-0.707)(0.897,-0.518)(0.966,-0.259)
\fi}
\def\DottedArc[#1]{
\ifcase#1
\bezier{4}(0.966,-0.259)(1.04,0)(0.966,0.259)
\or
\bezier{4}(0.966,0.259)(0.897,0.518)(0.707,0.707)
\or
\bezier{4}(0.707,0.707)(0.518,0.897)(0.259,0.966)
\or
\bezier{4}(0.259,0.966)(0,1.04)(-0.259,0.966)
\or
\bezier{4}(-0.259,0.966)(-0.518,0.897)(-0.707,0.707)
\or
\bezier{4}(-0.707,0.707)(-0.897,0.518)(-0.966,0.259)
\or
\bezier{4}(-0.966,0.259)(-1.04,0)(-0.966,-0.259)
\or
\bezier{4}(-0.966,-0.259)(-0.897,-0.518)(-0.707,-0.707)
\or
\bezier{4}(-0.707,-0.707)(-0.518,-0.897)(-0.259,-0.966)
\or
\bezier{4}(-0.259,-0.966)(0,-1.04)(0.259,-0.966)
\or
\bezier{4}(0.259,-0.966)(0.518,-0.897)(0.707,-0.707)
\or
\bezier{4}(0.707,-0.707)(0.897,-0.518)(0.966,-0.259)
\fi}
\def\Chord[#1,#2]{
\thinlines
\ifnum#1>#2\Chord[#2,#1]
\else\ifnum#1<#2
\ifcase#1
\ifcase#2
\or\qbezier(1,0)(0.516,0.138)(0.866,0.5)
\or\qbezier(1,0)(0.45,0.26)(0.5,0.866)
\or\qbezier(1,0)(0.327,0.327)(0,1)
\or\qbezier(1,0)(0.179,0.311)(-0.5,0.866)
\or\qbezier(1,0)(0.0536,0.2)(-0.866,0.5)
\or\put(1, 0){\line(-2, 0){2}}
\or\qbezier(1,0)(0.0536,-0.2)(-0.866,-0.5)
\or\qbezier(1,0)(0.179,-0.311)(-0.5,-0.866)
\or\qbezier(1,0)(0.327,-0.327)(0,-1)
\or\qbezier(1,0)(0.45,-0.26)(0.5,-0.866)
\or\qbezier(1,0)(0.516,-0.138)(0.866,-0.5)
\fi
\or\ifcase#2\or
\or\qbezier(0.866,0.5)(0.378,0.378)(0.5,0.866)
\or\qbezier(0.866,0.5)(0.26,0.45)(0,1)
\or\qbezier(0.866,0.5)(0.12,0.446)(-0.5,0.866)
\or\qbezier(0.866,0.5)(0,0.359)(-0.866,0.5)
\or\qbezier(0.866,0.5)(-0.0536,0.2)(-1,0)
\or\put(0.866, 0.5){\line(-5, -3){1.73}}
\or\qbezier(0.866,0.5)(0.146,-0.146)(-0.5,-0.866)
\or\qbezier(0.866,0.5)(0.311,-0.179)(0,-1)
\or\qbezier(0.866,0.5)(0.446,-0.12)(0.5,-0.866)
\or\qbezier(0.866,0.5)(0.52,0)(0.866,-0.5)
\fi
\or\ifcase#2\or\or
\or\qbezier(0.5,0.866)(0.138,0.516)(0,1)
\or\qbezier(0.5,0.866)(0,0.52)(-0.5,0.866)
\or\qbezier(0.5,0.866)(-0.12,0.446)(-0.866,0.5)
\or\qbezier(0.5,0.866)(-0.179,0.311)(-1,0)
\or\qbezier(0.5,0.866)(-0.146,0.146)(-0.866,-0.5)
\or\put(0.5, 0.866){\line(-3, -5){1}}
\or\qbezier(0.5,0.866)(0.2,-0.0536)(0,-1)
\or\qbezier(0.5,0.866)(0.359,0)(0.5,-0.866)
\or\qbezier(0.5,0.866)(0.446,0.12)(0.866,-0.5)
\fi
\or\ifcase#2\or\or\or
\or\qbezier(0,1.)(-0.138,0.516)(-0.5,0.866)
\or\qbezier(0,1.)(-0.26,0.45)(-0.866,0.5)
\or\qbezier(0,1.)(-0.327,0.327)(-1,0)
\or\qbezier(0,1.)(-0.311,0.179)(-0.866,-0.5)
\or\qbezier(0,1.)(-0.2,0.0536)(-0.5,-0.866)
\or\put(0, 1){\line(0, -2){2}}
\or\qbezier(0,1.)(0.2,0.0536)(0.5,-0.866)
\or\qbezier(0,1.)(0.311,0.179)(0.866,-0.5)
\fi
\or\ifcase#2\or\or\or\or
\or\qbezier(-0.5,0.866)(-0.378,0.378)(-0.866,0.5)
\or\qbezier(-0.5,0.866)(-0.45,0.26)(-1,0)
\or\qbezier(-0.5,0.866)(-0.446,0.12)(-0.866,-0.5)
\or\qbezier(-0.5,0.866)(-0.359,0)(-0.5,-0.866)
\or\qbezier(-0.5,0.866)(-0.2,-0.0536)(0,-1)
\or\put(-0.5, 0.866){\line(3, -5){1}}
\or\qbezier(-0.5,0.866)(0.146,0.146)(0.866,-0.5)
\fi
\or\ifcase#2\or\or\or\or\or
\or\qbezier(-0.866,0.5)(-0.516,0.138)(-1,0)
\or\qbezier(-0.866,0.5)(-0.52,0)(-0.866,-0.5)
\or\qbezier(-0.866,0.5)(-0.446,-0.12)(-0.5,-0.866)
\or\qbezier(-0.866,0.5)(-0.311,-0.179)(0,-1)
\or\qbezier(-0.866,0.5)(-0.146,-0.146)(0.5,-0.866)
\or\put(-0.866, 0.5){\line(5, -3){1.73}}
\fi
\or\ifcase#2\or\or\or\or\or\or
\or\qbezier(-1,0)(-0.516,-0.138)(-0.866,-0.5)
\or\qbezier(-1,0)(-0.45,-0.26)(-0.5,-0.866)
\or\qbezier(-1,0)(-0.327,-0.327)(0,-1)
\or\qbezier(-1,0)(-0.179,-0.311)(0.5,-0.866)
\or\qbezier(-1,0)(-0.0536,-0.2)(0.866,-0.5)
\fi
\or\ifcase#2\or\or\or\or\or\or\or
\or\qbezier(-0.866,-0.5)(-0.378,-0.378)(-0.5,-0.866)
\or\qbezier(-0.866,-0.5)(-0.26,-0.45)(0,-1)
\or\qbezier(-0.866,-0.5)(-0.12,-0.446)(0.5,-0.866)
\or\qbezier(-0.866,-0.5)(0,-0.359)(0.866,-0.5)
\fi
\or\ifcase#2\or\or\or\or\or\or\or\or
\or\qbezier(-0.5,-0.866)(-0.138,-0.516)(0,-1)
\or\qbezier(-0.5,-0.866)(0,-0.52)(0.5,-0.866)
\or\qbezier(-0.5,-0.866)(0.12,-0.446)(0.866,-0.5)
\fi
\or\ifcase#2\or\or\or\or\or\or\or\or\or
\or\qbezier(0,-1.)(0.138,-0.516)(0.5,-0.866)
\or\qbezier(0,-1.)(0.26,-0.45)(0.866,-0.5)
\fi
\or\ifcase#2\or\or\or\or\or\or\or\or\or\or
\or\qbezier(0.5,-0.866)(0.378,-0.378)(0.866,-0.5)
\fi\fi\fi\fi}
\def\FullChord[#1,#2]{
\Endpoint[#1]
\Endpoint[#2]
\Arc[#1]
\Arc[#2]
\Chord[#1,#2]
}
\def\EndChord[#1,#2]{
\Endpoint[#1]
\Endpoint[#2]
\Chord[#1,#2]
}
\def\Picture#1{
\begin{picture}(2,1)(-1,-0.167)
#1
\end{picture}
}
\def\DottedChordDiagram[#1,#2]{
\Picture{\DottedCircle \FullChord[#1,#2]}
}
\def\qq{\begin{eqnarray}}
\def\qqq{\end{eqnarray}}
\def\rx#1{~(\ref{#1})}
\def\ex#1{eq.\rx{#1}}
\def\eex#1{eqs.\rx{#1}}
\def\cx#1{~\cite{#1}}
\def\rw#1{~\ref{#1}}

\newlength{\shiftwidth}
\def\shift#1{&&\hbox to \shiftwidth{\hfill $\displaystyle#1$}}
\newlength{\sshiftwidth}
\def\sshift#1{\lefteqn{\hbox to
\sshiftwidth{\hfill$\displaystyle#1$}}}
\def\ie{{\it i.e.\ }}
\def\eg{{\it e.g.\ }}

\def\rhs{{\it r.h.s.\ }}
\def\lhs{{\it l.h.s.\ }}

\def\Tr{\mathop{{\rm Tr}}\nolimits}

\def\deg{ \mathop{{\rm deg}}\nolimits }
\def\p{^{\prime}}

\def\max{\mathop{{\rm max}}\nolimits}

\def\mmax#1{\max\{#1\}}

\def\Pexp{\mathop{{\rm Pexp}}\nolimits}
\def\PexpA#1{\Pexp \left(\oint_{#1} A_\mu dx^\mu \right)}

\def\pr#1#2{ \noindent{\em Proof of #1~\ref{#2}.} }

\def\qed{ \hfill $\Box$ }

\def\lrbc#1{ \left( #1 \right) }
\def\lrbs#1{ \left[ #1 \right] }

\def\spq{ \check{q} }
\def\xh{ h }

\def\suq{$SU_q(2)$ }

\def\va{ V_\a }
\def\vatn{ \va^{\otimes N} }
\def\vatno{ \va^{\otimes (N-1) } }
\def\vai{ V_{\a,\infty} }
\def\vaitn{ \vai^{\otimes N} }
\def\vaitno{ \vai^{\otimes (N-1)} }
\def\vava{ \va\otimes\va }
\def\vne{ V_N^{(\eta)} }
\def\vnoe{ V_{N-1}^{(\eta)} }
\def\vte{ V_2^{(\eta)} }
\def\vaivai{ \vai\otimes\vai }
\def\vno{ V_N^{(1)} }
\def\vnomo{ V^{(1)}_{N-1} }
\def\fovnomo{ f_0\otimes\vnomo }

\def\jakq{ J_\a(\cK;\spq) }
\def\vakq{ V_\a(\cK;\spq) }
\def\junk{ J_\a( {\rm unknot}; \spq ) }
\def\vnkz{ V^{(n)} (\cK;z) }
\def\vzkz{ V^{(0)}(\cK;z) }

\def\pnk#1{ P^{(n)}(\cK;#1) }
\def\pnkz{ \pnk{z} }

\def\aka#1{ \nabla_A(\cK;#1) }
\def\akuz{ \nabla_A({\rm unknot}; z) }
\def\akz{ \aka{z} }
\def\atnka#1{\nabla_A^{2n+1}(\cK;#1) }
\def\atnkz{ \atnka{z} }

\def\qd#1{ \spq^{#1\over 2} - \spq^{-{#1\over 2} } }
\def\qda{ \qd{\a} }
\def\qa{ \spq^{\a} }
\def\qma{ \spq^{-\a} }
\def\qfao{ \spq^{ {1\over 4} (\a^2-1)} }
\def\qfaom{ \spq^{ -{1\over 4} (\a^2-1)} }
\def\qmat{ \spq^{-{1\over 2}(\a-1) } }
\def\qat{ \spq^{ {1\over 2}(\a-1)} }
\def\qmaf{ \spq^{-{\a-1\over 4} } }
\def\qaf{ \spq^{ \a-1\over 4} }

\def\qhamn{ \spq^{ -{1\over 2}(\a-1)(m_1+m_2-n) } }
\def\qht{ \spq^{H\over 2} }
\def\qhttn{ \left( \qht \right)^{\otimes N} }
\def\qhttno{ \left( \qht \right)^{\otimes (N-1) } }
\def\qateb{ \spq^{ - {1\over 4}(\a^2-1)\,e(B_N)} }
\def\qatn{ \spq^{ {1\over 2} (\a-1)(N-1) } }
\def\qmet{ \spq^{-\eta} }
\def\qhane{ \spq^{ {1\over 2} (\a-1) (N-1-e(B_N)) } }
\def\qahebn{ \lrbc{\qa}^{ {1\over 2} (e(B_N) - N + 1) }  }
\def\qahneb{ \lrbc{\qa}^{ {1\over 2} (N - 1 - e(B_N)) }  }

\def\ohhne{ (1+h)^{-{1\over 2}(N-1-e(B_N))} }

\def\iqhtno{ I\otimes\qhttno }
\def\iqhtnobh{ \left( \iqhtno \right) \Bhn }

\def\zqqi{ \ZZ [ \spq, \spq^{-1} ] }
\def\zzt{ \ZZ[z^2] }
\def\zztt{ \ZZ[[z^2]] }
\def\Czn{ \IC[\zn] }
\def\Cztn{ \IC[\ztn] }
\def\Czt{ \IC[z_1,z_2] }
\def\Qmn{ \IQ[m,n] }
\def\zzlong{ \ZZ[\{a\}, \{a\p\}, \{e^{\e}\}, \qa, \qma] }
\def\Qqa{ \IQ[\qa,\qma] }
\def\Qdas{ \IQ[(\qda)^2] }
\def\zqqa{ \ZZ[\qa,\qma] }
\def\Zdas{ \ZZ[(\qda)^2] }
\def\zqqah{ \zqqa[[h]] }
\def\zh{ \ZZ[\xh] }
\def\zhh{ \ZZ[[\xh]] }

\def\oteoi{ \bigoplus_{\eta=0}^{\infty} }

\def\snz{ \sum_{n\geq 0} }
\def\szmn{ \sum_{0\leq m \leq n} }
\def\szmnt{ \sum_{0\leq m \leq {n\over 2} } }
\def\smz{ \sum_{m\geq 0} }
\def\scomb{ \sum_{0\leq n \leq M \atop 0 \leq m \leq {n\over 2} } }
\def\sznM{ \sum_{0\leq n \leq M} }
\def\sznmo{ \sum_{0\leq n \leq m_1} }
\def\setz{ \sum_{\eta\geq 0} }
\def\stjN{ \sum_{2\leq j\leq N} }
\def\setzmn{ \sum_{\etzmn} }
\def\setzmnmax{ \sum_{\etzmnmax} }
\def\sojM{ \sum_{1\leq j\leq M} }
\def\skz{ \sum_{k\geq 0} }
\def\sujo{ \sum_{j\geq 1} }
\def\setgmn{ \sum_{\etgmn} }
\def\sojNo{ \sum_{1\leq j\leq N-1} }
\def\szjN{ \sum_{0\leq j\leq N} }
\def\szktn{ \sum_{\zktn} }
\def\slon{ \sum_{1\leq l\leq n} }
\def\skpo{ \sum_{k\p \geq 1} }
\def\sko{ \sum_{k\geq 1} }
\def\sjz{ \sum_{j\geq 0} }
\def\smotl{ \sum_{m_1+1\leq l \leq m_2} }
\def\szkpk{ \sum_{0\leq k\p\leq k} }
\def\smmot{ \sum_{m_1\leq m\leq m_2} }
\def\setgmnmax{ \sum_{\etgmnmax} }

\def\plon{ \prod_{1\leq l \leq n} }
\def\plmtn{ \prod_{m_2+1\leq l \leq m_2+n} }
\def\plmoon{ \prod_{m_1+1\leq l \leq m_1+n} }
\def\plmon{ \prod_{m_1-n+1\leq l \leq m_1} }
\def\plmttn{ \prod_{m_2-n+1\leq l \leq m_2} }
\def\plzj{ \prod_{0\leq l\leq j-1} }
\def\plmnm{ \prod_{m-n+1\leq l \leq m} }
\def\pmzk{ \prod_{0\leq m \leq k} }
\def\plmpon{ \prod_{m+1\leq l\leq m+n} }
\def\poij{ \prod_{1\leq i \leq j} }
\def\pzijo{ \prod_{0\leq i\leq j-1} }
\def\pzik{ \prod_{0\leq i\leq k} }
\def\pzikkp{ \prod_{0\leq i\leq k-k\p} }
\def\pzij{ \prod_{0\leq i\leq j} }

\def\plnmc#1{ {\plmon \left( #1^l - 1 \right) \over \plon \left(
     #1^l - 1\right) } }
\def\plmn#1{ \plmtn \left(#1^{-l} - \spq^{-\a} \right) }

\def\plnmci#1{ {\plmttn \left( #1^{-l} - 1 \right) \over \plon \left(
     #1^{-l} - 1\right) } }
\def\plmni#1{ \plmoon \left(#1^{l} - \spq^{\a} \right) }
\def\plmonf{ {\plmon l\over n!} }
\def\plmtnf{ {\plmttn l\over n!} }
\def\plnmwc{ {\plmnm \lrbc{\ohlo } \over \plon \lrbc{\ohlo} } }
\def\plnmnf{ {\plmnm l\over n!} }

\def\ohlo{ (1+\xh)^l - 1 }
\def\ohlmo{ (1+\xh)^{-l} - 1 }

\def\ojN{ 1\leq j\leq N }
\def\tjN{ 2\leq j\leq N }
\def\etgNam{ \eta > (\a-1)(N-1) }
\def\etzmn{ 0\leq \eta \leq \nom}
\def\etzmnmax{ 0\leq \eta \leq (N-1)\amax}
\def\etgmn{ \eta > \nom }
\def\zetna{ 0 \leq \eta \leq (N-1)(\a-1) }
\def\etgmnmax{ \eta > (N-1) \amax }

\def\jgo{ j\geq 1 }
\def\mollmt{ m_1+1\leq l_1<\cdots<l_j\leq m_2 }

\def\Dmnk{ \Dmn (\cK) }
\def\dmnk{ \dmn (\cK) }
\def\Dmn{ D_{m,n} }
\def\dmn{ d^{(n)}_m }

\def\ohat{ (1+\xh)^{\a\over 2} - (1+\xh)^{-{\a\over 2}} }

\def\Rchk{ \check{R} }
\def\Rpchk{ \check{\cR} }
\def\Rpchki{ \Rpchk^{({\rm inv})} }
\def\Rpchke{ \Rpchk [a;\e_1,\e_2,\e_{12}] }
\def\Rpchkie{ \Rpchki [a\p;\e_1,\e_2,\e_{12}] }
\def\Rhpchk{ \tilde{ \check{\cR} } }
\def\Rhpchke{ \Rhpchk[a,\e_1,\e_2,\e_{12} ] }
\def\Rhpchki{ \tilde{ \check{\cR} }^{({\rm inv})} }
\def\Rhpchkie{ \Rhpchki[a\p,\e_1,\e_2,\e_{12} ] }

\def\Bhn{ \hat{B}_N }
\def\Bhno{ \Bhn^{(1)} }
\def\Bhnp{ \hat{\cB}_N }
\def\Bhnpe{ \Bhnp[\spar] }
\def\Bhtnp{ \tilde{ \hat{\cB} }_N }
\def\Bhtnpe{ \Bhtnp[\spar] }

\def\Tje{ T_j(B_N | \spar) }
\def\Tjed{ T_j(B_N | \spard) }
\def\Tjkro{ T_{j,k}^{(R,1)} }
\def\Tjokro{ T_{j+1,k}^{(R,1)} }
\def\Tjkpro{ T_{j,k\p}^{(R,1)} }
\def\Tjrt{ T_j^{(R,2)} }
\def\tTjrt{ \tilde{T}_j^{(R,2)} }
\def\Tojk{ T^{(1)}_{j+1,k,k\p} }
\def\tTojk{ \tilde{T}^{(1)}_{j+1,k,k\p} }
\def\Trro{ T^{(R,1)} }
\def\Trt{ T^{(R,2)} }

\def\tebn{ T_\eta(B_N) }
\def\ten{ T_{\eta,n} }
\def\tene{ \ten(B_N|\qa,\qma) }
\def\vnee{ V_n(B_N|\qa,\qma,\lambda) }

\def\tQj{ \tilde{Q}_j }
\def\tQje{ \tQj(\spar) }
\def\qnk{ Q_{n,k} }
\def\qnke{ \qnk(\qa,\qma) }
\def\qje{ Q_j(\qa,\qma) }

\def\Bnm{ B_n(m_1,m_2) }

\def\pnqqa{ P_n(\qa,\qma) }

\def\Tro{ \Tr^{(1)} }
\def\Trfov{ \Tr_{\fov} }
\def\Trfovi{ \Tr_{\fovi} }
\def\Trfovnoe{ \Tr_{\fovnoe} }
\def\Trvne{ \Tr_{\vne} }

\def\detvno{ \det_{\vno} }
\def\detfovnomo{ \det\nolimits_{\fovnomo} }
\def\detfovnoe{ \det\nolimits_{\fovnoe} }

\def\fmot{ f_{m_1}\otimes f_{m_2} }
\def\fmoti{ f_{m_2}\otimes f_{m_1} }
\def\fmtn{ f_{m_1}\otimes \cdots \otimes f_{m_N} }
\def\fmnot{ f_{m_1-n} \otimes f_{m_2+n} }
\def\fmnoti{ f_{m_2+n} \otimes f_{m_1-n} }
\def\fmnotii{ f_{m_2-n} \otimes f_{m_1+n} }
\def\fov{ f_0\otimes \vatno }
\def\fovi{ f_0\otimes \vaitno }
\def\fovnoe{ f_0\otimes \vnoe }
\def\fzmtn{ f_0\otimes f_{m_2} \otimes \cdots \otimes f_{m_N} }
\def\fovcap{ \fov\bigcap \fovnoe }

\def\zmot{ z_1^{m_1}z_2^{m_2} }
\def\zmnoti{ z_1^{m_2+n} z_2^{m_1-n} }

\def\zmn{ z_1^{m_1} \cdots z_N^{m_N} }
\def\zn{ z_1, \ldots, z_N }
\def\ztn{ z_2, \ldots, z_N }

\def\eon{ e_1,\ldots,e_N }

\def\zmao{ 0\leq m\leq \a-1 }
\def\zmtao{ 0\leq m_2\leq \a-1 }
\def\mga{ m\geq \a }
\def\mtna{ m_2+n\geq \a }
\def\zktn{ 0\leq k \leq 2n(N-1) }

\def\ffd{ \lrbc{ 1 + \sojM \hj \Tjed } }
\def\fsd{ \lrbc{ 1 + \sojM \hj {(-1)^j\over j!} \partial_\kappa^j } }

\def\bckae{ \right|_{
  {{{\kappa=1 \hfill\atop \{a\}=1-\qma\hfill}\hfill\atop \{a\p\}=1-\qa
   \hfill} \hfill\atop \{\e\}=0 \hfill} } }
\def\bcaae{ \right|_{
  {{ \{a\}=1\hfill\atop \{a\p\}=1
   \hfill} \hfill\atop \{\e\}=0 \hfill} } }
\def\bcae{ \right|_{a=1-\qma\hfill \atop \e_1=\e_2=\e_{12}=0} }
\def\bclz{ \right|_{\lambda=0} }
\def\bcllz{ \right|_{\lambda=\lambda_0} }
\def\bcaaz{ \right|_{\{a\}=0 \hfill\atop \{a\p\}=0\hfill} }
\def\bcaez{ \right|_{a=1-\qma\hfill\atop \e_1=\e_2=\e_{12}=0\hfill} }
\def\bcapez{ \right|_{a\p=1-\qa\hfill\atop
     \e_1=\e_2=\e_{12}=0\hfill} }
\def\bcaape{ \right|_{
  {{\{a\}=1-\qma\hfill\atop \{a\p\}=1-\qa
   \hfill} \hfill\atop \{\e\}=0 \hfill} } }

\def\fone#1#2#3{1+
\sujo\xh^j\,
 {\plzj(#1-l)\over (1-\qma)^j}
     \skz \xh^k\,\Tjkro(#2,#3) }
\def\ftwo#1#2{1+\sujo\xh^j\,\Tjrt(#1,#2) }
\def\fthree#1{1+\sujo\xh^j\,{\plzj(#1-l)\over j!} }

\def\sjo{ \sigma_{j,j+1} }
\def\hnah{ {1\over 2} [N(\a-1) - H] }
\def\hah{ {1\over 2}(\a-1-H) }

\def\ohmo{ \cO(\xh^{M+1}) }
\def\nom{ (N-1)M }
\def\htz{ \xh\rightarrow 0 }
\def\poh{ (1+\xh) }
\def\spar{ \{a\},\{a\p\},\{\e\} }
\def\spard{ \{\partial_a\},\{\partial_{a\p} \}, \{\partial_\e \} }
\def\keta{ \kappa^\eta }
\def\leta{ \lambda^\eta }
\def\hj{ \xh^j }
\def\tcO{ \tilde{\cO} }
\def\fpar{ {1\over \eta!}\,\partial_\lambda^\eta }
\def\amax{ \mmax{ \a-1,M } }

\def\nabla{ \Delta }

\begin{document}
%\draft

\begin{titlepage}
%\centerline{\hfill                 UMTG-183-95}
\centerline{\hfill                 q-alg/9604005}
\vfill
\begin{center}

{\large \bf
The Universal $R$-Matrix, Burau Representation and the Melvin-Morton
Expansion of the Colored Jones Polynomial.
}
\\

\bigskip
\centerline{L. Rozansky
%\footnote{Work supported by the National Science Foundation
%under Grant No. PHY-92 09978.}
}

\centerline{\em School of Mathematics, Institute for Advanced Study}
\centerline{\em Princeton, NJ 08540, U.S.A.}
\centerline{{\em E-mail address: rozansky@math.ias.edu}}

\vfill
{\bf Abstract}

\end{center}
\begin{quotation}
P.~Melvin and H.~Morton\cx{MeMo} studied the expansion of the colored
Jones polynomial of a knot in powers of $\spq - 1$ and color. They
conjectured an upper bound on the power of color versus the power of
$\spq - 1$. They also conjectured that the bounding line in their
expansion generated the inverse Alexander-Conway polynomial. These
conjectures were proved by D.~Bar-Natan and S.~Garoufalidis\cx{BG}.

We have conjectured\cx{Ro7} that other `lines' in the Melvin-Morton
expansion are generated by rational functions with integer
coefficients whose denominators are powers of the Alexander-Conway
polynomial. Here we prove this conjecture by using the $R$-matrix
formula for the colored Jones polynomial and presenting the universal
$R$-matrix as a `perturbed' Burau matrix.
\end{quotation}
\vfill
\end{titlepage}

\pagebreak

%+++++++++++++++++++++++++++++++
\nsection{Introduction}
\label{s1}
%+++++++++++++++++++++++++++++++
\hyphenation{Re-she-ti-khin}
\hyphenation{Tu-ra-ev}

Let $\cK$ be a knot in $S^3$ endowed with canonical framing (\ie, its
self-linking number is zero). We assign to this knot an
$\a$-dimensional \suq module $\va$. $\jakq$ denotes the colored Jones
polynomial of $\cK$, normalized in such a way that it is
multiplicative under a disconnected sum and
\qq
\junk = {\qd{\a} \over \qd{1} }.
\label{1.1}
\qqq
Another popular normalization for the Jones polynomial is
\qq
\vakq = {\jakq \over \junk} = {\qd{1} \over \qd{\a} }\, \jakq.
\label{1.2}
\qqq
The advantage of the normalization\rx{1.2} is that
\qq
\vakq \in \zqqi.
\label{1.3}
\qqq

The colored Jones polynomial is an effective invariant of knots.
However, its relation to classical topological invariants of knots
remains mostly obscure (see, \eg a review\cx{Bi}). One may try to
decompose $\vakq$ into some simpler `building blocks' in a hope that
their topological nature would be easier to establish. An important
step in this direction was made by P.~Melvin and H.~Morton\cx{MeMo}.
They suggested to expand the Jones polynomial in powers of $\a$ and
\qq
\xh = \spq - 1
\label{1.4}
\qqq
(actually, they expanded $\jakq$ and used $\log(1+\xh)$ rather than
$\xh$ as an expansion parameter; this seemed to be a more `physically
natural' choice). If we expand $\vakq$ in Taylor series in $\xh$ for
a fixed value of $\a$ then, according to\cx{MeMo}, the coefficients
will be finite degree polynomials in $\a$:
\qq
\vakq = \snz \xh^n \left( \szmn \Dmnk \a^{2m} \right).
\label{1.5}
\qqq
Melvin and Morton proved that the coefficients $\Dmnk$ are finite
type invariants of $\cK$ of order $n$. Moreover, they conjectured
that the bound on powers of polynomials of $\a$ can be improved and
that the bounding line contribution to $\vakq$ is equal to the
inverse Alexander-Conway polynomial:
%%%%%%%%%%%
\begin{theorem}
\label{t1.1}
For a knot $\cK\subset S^3$ the coefficients $\Dmnk$ of the
expansion\rx{1.5} satisfy the following two properties:
\qq
&
\Dmnk = 0 \qquad \mbox{{\rm for} $ m > {n\over 2}$,}
\label{1.7}\\
&
\smz D_{m,2m}(\cK) \, a^{2m} =
{1\over \aka{e^{i \pi a} - e^{-i\pi a} } },
\label{1.8}
\qqq
here $a $ is a formal parameter and $\aka{z}$ is the Alexander-Conway
polynomial of $\cK$ which satisfies the skein relation of
Fig.~\ref{fig:skein}
and is normalized by the condition
\qq
\akuz = 1.
\label{1.9}
\qqq
\end{theorem}
%%%%%%%%%%%%%%%%%%%%%%%
This theorem was proved by D.~Bar-Natan and S.~Garoufalidis\cx{BG}.

\begin{figure}[hbt]
\def\KZero{
\begin{picture}(2,2)(-1,-1)
\put(0,0){\circle{2}}
\put(-0.707,-0.707){\vector(1,1){1.414}}
\put(0.707,-0.707){\vector(-1,1){1.414}}
\put(0,0){\circle*{0.15}}
\end{picture}}

\def\KPlus{ \Picture {
\put(0,0){\circle{2}}
\put(-0.707,-0.707){\vector(1,1){1.414}}
\put(0.707,-0.707){\line(-1,1){0.6}}
\put(-0.107,0.107){\vector(-1,1){0.6}}
} }
%\end{picture}}

%\def\KMinus{
%\begin{picture}(2,2)(-1,-1)

\def\KMinus{ \Picture {
\put(0,0){\circle{2}}
\put(-0.707,-0.707){\line(1,1){0.6}}
\put(0.107,0.107){\vector(1,1){0.6}}
\put(0.707,-0.707){\vector(-1,1){1.414}}
} }
%\end{picture}}

%\def\KII{
%\begin{picture}(2,2)(-1,-1)

\def\KII{ \Picture {
\put(0,0){\circle{2}}
\qbezier(-0.707,-0.707)(0,0)(-0.707,0.707)
\qbezier(0.707,-0.707)(0,0)(0.707,0.707)
\put(-0.607,0.607){\vector(-1,1){0.1414}}
\put(0.607,0.607){\vector(1,1){0.1414}}
} }
%\end{picture}}

\begin{displaymath}
{\nabla_A\left(\KPlus\; ;z\right)}~\qquad - \qquad
\nabla_A\left(\KMinus\; ;z\right)~\qquad = \qquad z\;
\nabla_A\left(\KII\; ;z\right)~.
\end{displaymath}
\caption{The skein relation for the Alexander-Conway polynomial}
\label{fig:skein}
\end{figure}

In\cx{Ro7} we conjectured that the expansion\rx{1.5} satisfies some
further properties. To formulate our conjecture we rearrange \ex{1.5}
as a formal power series in $\xh$ and $\a\xh$:
\qq
\vakq = \snz \xh^n \smz D_{m,n+2m}\,(\a\xh)^{2m}.
\label{1.10}
\qqq
Then, motivated by \ex{1.8}, we introduce a new variable
\qq
z = \qda
\label{1.11}
\qqq
instead of $\a\xh$:
\qq
\a\xh = 2\log\left( \sqrt{ 1 + \left( {z\over 2} \right)^2 }
+ {z\over 2} \right) {\xh \over \log(1+\xh)} =
z + \cO(z^3,\xh).
\label{1.12}
\qqq
Substituting \ex{1.12} into \ex{1.10} as a formal power series we get
the expansion in powers of $\xh$ and $z^2$:
\qq
&\vakq = \snz \vnkz \, \xh^n,
\label{1.13}\\
&\vnkz = \smz \dmnk \, z^{2m}.
\label{1.14}
\qqq

Eq.\rx{1.8} implies that the formal power series $\vzkz$ comes from
the Taylor expansion of the inverse Alexander-Conway polynomial at
$z=0$:
\qq
\vzkz = {1\over \akz}.
\label{1.15}
\qqq
We conjectured (`strong conjecture' of\cx{Ro7}) that the `upper
lines' $\vnkz$ of the expansion\rx{1.13} come also from the Taylor
expansion of rational functions of $z$.
%%%%%%%%%%%%%%%%%%%
\begin{theorem}
\label{t1.2}
There exists a set of polynomial
invariants of knots in $S^3$
\qq
\pnkz \in \zzt,
\label{1.16}
\qqq
such that the series $\vnkz$ is the Taylor series
expansion at $z=0$ of the rational function
\qq
\vnkz = {\pnkz \over \atnkz}.
\label{1.17}
\qqq
\end{theorem}
In\cx{Ro7} we presented experimental evidence in support of this
theorem by calculating the first polynomials $\pnkz$ of some simple
knots.

The Alexander-Conway polynomial of a knot satisfies the properties
\qq
&\akz \in \zzt,
\label{1.18}\\
& \aka{0} = 1.
\label{1.19}
\qqq
Therefore as a Taylor series expansion at $z=0$,
\qq
{1\over \akz} \in \zztt,
\label{1.20}
\qqq
and in view of \ex{1.16} we have the following simple corollary
 of the Theorem\rw{t1.2} (`weak conjecture' of\cx{Ro7}):
%%%%%%%%%%%%%%%%%
\begin{corollary}
\label{c1.1}
All the coefficients $\dmnk$ of the expansion\rx{1.14} are integer:
\qq
\dmnk \in \ZZ.
\qqq
\end{corollary}
%%%%%%%%%%%%%%%%%%
Note that this corollary is stronger than the condition
\qq
n!\,\Dmnk \in \ZZ,
\label{1.22}
\qqq
which comes from the fact that because of\rx{1.3}
\qq
\szmnt \Dmnk\, \a^{2m} \in \ZZ
\qquad \mbox{ for $\a\in\ZZ$.}
\label{1.23}
\qqq

R.~Lawrence\cx{Lw} has formulated a conjecture about $p$-adic
properties of Ohtsuki's invariants of rational homology spheres. We
proved that conjecture in\cx{Ro8} for the special case of a manifold
constructed by a rational surgery on a knot in $S^3$. Our proof was
based on Corollary\rw{c1.1} which we used as a conjecture. Thus the
proof of Theorem\rw{t1.1} that will be presented in this paper,
completes the proof of\cx{Ro8}.

We will prove Theorem\rw{t1.2} in the following equivalent form:
%%%%%%%%%%%%%%%%%%
\begin{proposition}
\label{p1.1}
For a knot $\cK\subset S^3$ there exists a set of polynomials
$\pnkz \in \zzt$ such that for a fixed $\a$ and for any $M>0$
\qq
\scomb \Dmnk\,\a^{2m}\xh^n =
\sznM \xh^n\, {\pnk{\ohat} \over \atnka{\ohat} }
+ \cO(\xh^{M+1}).
\label{1.023}
\qqq
\end{proposition}
The equivalence between this proposition and Theorem\rw{t1.2} follows
from the structure of the substitution\rx{1.12}. A coefficient
$\dmn$ of the expansion\rx{1.14} is determined only by the
coefficients $D_{m\p,n\p}$, $m\p\leq m$, $n\p\leq n+2m$. Therefore
\ex{1.023} for $M=n+2m$ indicates that the coefficient $\dmn$ does
come from the Taylor expansion of the \rhs of \ex{1.17} at $z=0$.

\noindent{\em Outline of the Proof of Proposition~\ref{p1.1} }.
Similarly to\cx{MeMo}, we use the expression for the colored Jones
polynomial $\jakq$ as a (quantum) trace of the product of
$\Rchk$-matrices. The matrices correspond to elementary braids in the
picture of the knot $\cK$ as a closure of an $N$-strand braid $B_N$.
The trace is taken over the tensor product $\vatn$ of the
$\a$-dimensional \suq modules $\va$. A choice of a particular
basis $f_m$, $\zmao$ in $\va$ allows us to map $\vatn$ into the
algebra $\Czn$ of polynomials of $N$ variables:
\qq
\fmtn \rightarrow \zmn.
\label{1.24}
\qqq
The action of $\Rchk$-matrices on $\vatn$ can be extended to the
space $\Czn$. It turns out that in the approximation, when
$\xh\rightarrow 0$ and $\qa$ is kept constant, the action
of the $\Rchk$-matrix becomes an endomorphism of the algebra $\Czn$.
This endomorphism is
%In other words, the action of the $\Rchk$-matrix on polynomials is
generated by a linear transformation of the space $\IC^N$ of
variables $\zn$. This transformation coincides with the Burau
representation of elementary braids.

Let $\cO$ be such
an endomorphism of $\Czn$. Its trace over the whole
algebra $\Czn$ can be expressed in terms of the
restriction $\tilde{\cO}$ of the action of $\cO$ to the subspace
$\IC^N$ of variables $\zn$. More precisely,
\qq
\Tr_{\Czn} \cO = {1\over \det_{\IC^N} (1-\cO)}.
\label{1.25}
\qqq
Since the matrix $\cO$ coming from the braid $B_N$ coincides with the
Burau representation, the determinant in the denominator of the \rhs
of \ex{1.25} is equal to the Alexander-Conway polynomial of $\cK$.
This leads to the Melvin-Morton relation\rx{1.15} Since we never use
the bound\rx{1.7}, these considerations constitute yet another,
$R$-matrix based, proof of the Melvin-Morton conjecture.

The exact $\Rchk$-matrix presented as a series in $\xh$ is not an
endomorphism of $\Czn$. The coefficients at higher powers of $\xh$
are rather `perturbed endomorphisms' whose traces can be calculated
with the help of standard tricks of Quantum Field Theory. These
tricks are rigorous in our finite-dimensional context. The result of
such calculation is the expansion\rx{1.13},\rx{1.14},\rx{1.17}.

The fact that Burau representation generates the $R$-matrix based
braid group representation in the limit of $\xh \rightarrow 0$,
$\qa = \mbox{const}$,
can be traced back to the equivalence between the $R$-matrix
representation and the action of the braid group on the twisted $m$th
cohomology of the configuration space of $m$ points on a complex
plane with $N$ holes, considered by R.~Lawrence\cx{Lw1},\cx{Lw2}
(see also the book\cx{Va},
this construction is known in physical literature as `free field
representation').

Although we try deliberately to avoid any `physical' references in
this paper, we can not help but mention that our calculations are
very similar to those of\cx{KS}. L.~Kauffman and H.~Saleur related
the Burau representation to the evolution of fermionic particles
along the strands of the braids. The particles scatter at the
elementary braids. The scattering is described by the Burau matrix.
It is essential that the particles are
free, that is, they scatter independently of one another. It follows
from our calculations that in the limit of $\xh\rightarrow 0$,
$\qa$-fixed, the colored Jones polynomial comes from
the evolution of `almost free' bosonic particles. The
`particles' are monomials $z_j$, $\ojN$ of $\Czn$, and their
`freedom' is the physical equivalent of the mathematical notion of
endomorphism. The leading term\rx{1.15} corresponds to absolutely
free particles and expresses the fact that fermions are `inverse'
bosons. The subsequent terms\rx{1.17} come from the perturbative
expansion in the weak coupling constant $\xh$.

Here is the plan of the paper.
In Section\rw{s2} we recall the formula for the colored Jones
polynomial as a trace of the product of $\Rchk$-matrices over the
tensor product $\vatn$ of \suq modules $\va$. We break a strand in
the braid closure and also extend the trace to the product of Verma
modules $\vaitn$. In Section\rw{s3} we discuss the expansion of the
$\Rchk$-matrix in the limit of $\xh\rightarrow 0$. We express the
terms of this expansion as derivatives of the `parametrized'
$\Rchk$-matrix. In Section\rw{s4} we present parametrized
$\Rchk$-matrix as an endomorphism of the algebra $\Czn$ and
calculate its trace. In Section\rw{s5} we combine the results of all
the previous sections and complete the proof of the
Proposition\rw{p1.1}. Discussion contains some comments about
generalizing our calculations to links. In Appendix we prove
some technical lemmas which were needed in
Section\rw{s3}.

%In Appendix\rw{A2} we comment briefly on the
%relation between $R$-matrix and Burau representations in the context
%of the work of R.~Lawrence\rw{Lw1},\rw{Lw2}.

%%%%%%%%%%%%%%%%%%%
\nsection{The colored Jones polynomial as a character of the braid
group}
\label{s2}

In this section we review the representation of the braid group based
on the universal $R$-matrix (\cx{RT1}, see also a nice exposition
in\cx{KM}). We also recall a formula for the colored Jones polynomial
as a trace of this braid representation with one broken closure
strand.

\subsection{The \suq $R$-matrix}
\label{ss2.1}
Let $\va$ be the $\a$-dimensional \suq module. We choose the basis
vectors $f_m$, $\zmao$ of $\va$ in such a way that the action of the
standard \suq generators is
\qq
X\, f_m & = & [m]\,f_{m-1},
\label{2.1}\\
Y\, f_m & = & [\a-1-m]\, f_{m+1},
\label{2.2}\\
H\, f_m & = & (\a-1-2m) \, f_m,
\label{2.3}
\qqq
here we use the notation
\qq
[n] = {\qd{n}\over \qd{1} }.
\label{2.4}
\qqq
Our basis vectors $f_m$ are related to the basis vectors $e_j$
of\cx{KM}, eq.(2.8):
\qq
f_m = { [\a-m-1]!\, [m]! \over [\a-1]!}\; e_{ {\a-1\over 2} -m},
\label{2.5}
\qqq
here
\qq
[n]! =
\left\{
\begin{array}{cl}
\plon [l] &\mbox{if $n\geq 1$}\\
1&\mbox{if $n=0$}.
\end{array}
\right.
\label{2.6}
\qqq
The action of the $R$-matrix on the basis vectors $\fmot$ of the
tensor product $\va\otimes\va$ is given by the formula which comes
from substituting \ex{2.5} into Corollary~2.3.2 of\cx{KM}:
\qq
R(\fmot) & = &
\sznmo \left( \qd{1} \right)^n
{ [\a - m_2 - 1]! \over [\a-m_2-n-1]! }
\,{[m_1]! \over [m_1-n]! \, [n]!}
\label{2.7}\\
&&\qquad
\times
\spq^{{1\over 4} [\a-2m_1-1)(\a-2m_2-1) + 2n(m_1-m_2) - n(n+1)]}
\,\fmnot.
\nonumber
\qqq
After some simple transformations this formula becomes
\qq
R(\fmot) & = &
\qfao \qmat \snz \;\plmn{\spq} \plnmc{\spq}
\nonumber
\\
&&\qquad\times
\qhamn \,\fmnot.
\label{2.8}
\qqq
Note that in our conventions $\sum_{a\leq j \leq b} (\cdots) = 0$
and $\prod_{a\leq j\leq b} (\cdots) = 1$ if $a>b$.

The `flipped' matrix $\Rchk$ is defined as
\qq
\Rchk = PR,
\label{2.08}
\qqq
here $P$ is the permutation operator:
\qq
P\,(\fmot) = \fmoti.
\label{2.9}
\qqq
In order to simplify the future calculations we will write the matrix
elements of $\Rchk$ in the basis with `rotated phases':
\qq
\fmot \rightarrow \left( \qat \right)^{m_2} \fmot.
\label{2.10}
\qqq
In other words, our $\Rchk$-matrix is defined as
\qq
\Rchk = \left( \qmaf \right)^{I\otimes H} PR\,
\left(\qaf\right)^{I\otimes H}
\label{2.11}
\qqq
($I$ is the identity operator) instead of \ex{2.8}:
\qq
\Rchk(\fmot) & = &
\qfao \qmat \snz \; \plmn{\spq} \plnmc{\spq}
\nonumber\\
&&\qquad\times
\spq^{-(\a-1)m_2} \,\fmnoti.
\label{2.12}
\qqq
The elements of the inverse $\Rchk$-matrix can be obtained with the
help of relation
\qq
\Rchk^{-1}(\spq) = P\,\Rchk(\spq^{-1})\,P.
\label{2.012}
\qqq
In other words, we have to substitute
\qq
\fmot \rightarrow \fmoti,\qquad
\fmnoti \rightarrow \fmnot
\label{2.1012}
\qqq
and $\spq \rightarrow \spq^{-1}$ in \ex{2.12}. As a result,
\qq
\Rchk^{-1}(\fmot) & = &
\qfaom \qat \snz \; \plmni{\spq} \plnmci{\spq}
\nonumber\\
&&\qquad\times
\spq^{(\a-1)m_1}\,\fmnotii.
\label{2.13}
\qqq
\subsection{Braid group representation and its character}
\label{ss2.2}

Let $B_N$ be a braid of $N$ strands. We associate with it a tensor
product $\vatn$ of $N$ \suq modules $\va$, one module per each
position. To an elementary positive braid $\sjo$ which switches the
strands at $j$th and $(j+1)$st positions we associate the
$\Rchk$-matrix acting on $\vava$ at $j$th and $(j+1)$st positions in
$\vatn$.
The operator
$\Rchk^{-1}$ is associated to a negative elementary braid
$\sigma^{-1}$. This construction is known to represent the braid
group in $\vatn$ in such a way that the action of braids commutes
with the `global' (\ie, acting simultaneously on all individual
spaces $\va$ in $\vatn$) action of \suq. We denote the representation
of a brain $B_N$ in $\vatn$ as $\Bhn$. Note that because of
\ex{2.11}, our representation is a conjugation of the standard one by
the operator
\qq
\left(\qaf\right)^{H} \otimes
\left(\qaf\right)^{2H} \otimes
\cdots \otimes
\left(\qaf\right)^{NH}.
\label{2.013}
\qqq

Suppose that a knot $\cK\subset S^3$ is presented as a closure of a
braid $B_N$. The colored Jones polynomial of $\cK$ can be calculated
as a `quatum trace' of $\Bhn$. We choose the closure strands to go to
the right of the braid. Denote by $\qhttn$ an operator that acts as
$\qht$ on every module $\va$ of $\vatn$. Then
\qq
\jakq = \qateb \Tr_{\vatn} \qhttn \Bhn,
\label{2.14}
\qqq
here $e(B_N)$ is the number of positive elementary braids minus the
number of negative elementary braids in $B_N$. The prefactor $\qateb$
is the framing correction. It is due to the fact that the knot $\cK$
constructed by closing the braid $B_N$, has the blackboard framing
$e(B_N)$.

\subsection{Breaking a closure strand}
\label{ss2.3}

We are going to use the $SU(2)_q$-invariance of the braid group
representation in order to reduce the trace of \ex{2.14}.

Consider a tangle constructed by closing all braid positions except
the first one which we leave open. To this tangle we associate an
operator $\Bhno$ acting on $\va$ (the first element in the product
$\vatn$):
\qq
\Bhno = \Tro_{\vatno} \iqhtnobh.
\label{2.15}
\qqq
The operator $\iqhtno$ in this formula acts as identity on
the first $\va$ and as $\qht$ on all other $\va$ of $\vatn$.
$\Tr_{\vatno}$ is the trace taken over all $\va$ of $\vatn$ except
the first one.

The \suq invariance of $\Bhno$ together with irreducibility of $\va$
means that $\Bhno$ is proportional to the identity operator:
\qq
\Bhno = CI, \qquad C\in \IC.
\label{2.16}
\qqq
On the other hand, in view of \eex{2.14} and\rx{2.15}, the Jones
polynomial $\jakq$ is equal to the quantum trace of $\Bhno$ which
corresponds to closing the remaining strand:
\qq
\jakq & = &
\qateb \Tr_{\va} \qht \Bhno
= \qateb \,C \Tr_{\va} \qht
\nonumber\\
& = &
\qateb\, C\,{\qd{\a}\over \qd{1} }.
\label{2.17}
\qqq
Then, according to \eex{1.1},\rx{1.2},
\qq
\vakq = \qateb\,C.
\label{2.017}
\qqq

To find the constant $C$ we can choose any diagonal matrix element of
$\Bhno$. We will do it for $f_0$. Thus
\qq
\vakq = \qateb \,\Trfov \iqhtnobh.
\label{2.18}
\qqq
The symbol $\Trfov$ has the following meaning. The operator
\qq
\iqhtnobh
\label{2.018}
\qqq
acts on the full space $\vatn$. We project this action
onto the subspace $\fov \subset \vatn$ along the other subspaces
$f_m\otimes\vatno$, $m\geq 1$. Then we take the trace of this
projection. In other words, we simply take the diagonal matrix
element of\rx{2.018} for $f_0$ (with respect to our basis $f_m$) in
the first space $\va$ and take traces over all other spaces $\va$ of
$\vatn$.

\subsection{Extension and stratification of the trace}
\label{ss2.4}

Our next step is to extend the \suq modules $\va$ to the
infinite-dimensional spaces $\vai$ by adding formally the basis
vectors $f_m$, $\mga$ to the already existing vectors $f_m$, $\zmao$.
The matrix elements of \suq generators $X,Y,H$ are still defined by
\eex{2.1},\rx{2.2} and\rx{2.3}, while the action of $\Rchk$-matrix is
given by \ex{2.12}.

The tensor product $\vaitn$ can be decomposed into a direct sum of
linear spaces $\vne$ corresponding to the eigenvalues $\eta$ of the
operator $\hnah$ acting on $\vaitn$:
\qq
\vaitn = \oteoi \vne.
\label{2.21}
\qqq
The subspace $\fovi$ is also decomposed:
\qq
\fovi = \oteoi \fovnoe, \qquad \fovnoe \subset \vne.
\label{2.22}
\qqq
The spaces $\vne$, $\fovnoe$ are finite-dimensional.

Since $H$ commutes with $\Rchk$-matrices\rx{2.12} and with the
operator $\iqhtno$, we can split the trace\rx{2.18} into the traces
over the eigenspaces of $H$.
%%%%%%%%%%%%%%%%%%%%%%%%%%%%
\begin{proposition}
\label{p2.1}
For any $M>0$,
\qq
\vakq
=\qateb \,\qatn
%\!\!\!\!\!\!\!\!
\setzmnmax
%\!\!\!\!
\qmet \Trfovnoe \Bhn.
%\nonumber\\
%&&\qquad\qquad\qquad\qquad\qquad
%\nonumber
\label{2.26}
\qqq
\end{proposition}
%%%%%%%%%%%%%%%%%%%
\pr{Proposition}{p2.1}
First of all, since
\qq
\fov \subset \bigoplus_{\zetna} \fovnoe
\label{x.1}
\qqq
and since $\fovnoe$ is an eigenspace of the operator $\iqhtno$ with
the eigenvalue $\qatn \qmet$, we can transform \ex{2.18} into
\qq
\vakq
=\qateb \,\qatn
\!\!\!\!\!\!\!\!
\setzmnmax
\!\!\!\!
\qmet \Tr_{\fovcap} \Bhn.
\label{x.2}
\qqq

It remains to check that
\qq
\Tr_{\fovcap} \Bhn = \Trfovnoe \Bhn.
\label{x.3}
\qqq
To show this consider
a matrix element of the $\Rchk$-matrix between the vectors
$\fmot$ and $\fmnoti$ such that $\zmtao$ and $\mtna$. In this case
the product
\qq
\plmn{\spq}
\label{2.20}
\qqq
of \ex{2.12} contains a term corresponding to $l=\a$, so the matrix
element is zero. Therefore if we follow the evolution of a vector
{\em along the braid strand}, then we see that transitions from the
elements $f_m$, $\zmao$ to $f_m$, $\mga$ are forbidden at positive
elementary braids. The same analysis of \ex{2.13} shows that these
transitions are also forbidden at negative braids.

In the calculation of $\Trfovnoe$ the evolution of the vector along
the braid strands and the closure strands starts at $f_0$ at the
beginning of the first strand. Following this evolution we will
cover all the segments of the braid, because its closure is a knot.
Therefore the elements $f_m$, $\mga$ will never appear in the
calculation of the \rhs of \ex{x.3}. \qed

%%%%%%%%%%%%%%%%%%%%%%%%%%%%
\nsection{Expansion of $\Rchk$-matrix}
\label{s3}

In order to use \ex{2.26} for calculation of the coefficients $\Dmn$
we have to expand the matrix elements of $\Bhn$, which appear in the
\rhs of \ex{2.26}, in powers of $\xh$ at $\xh=0$. These matrix
elements come from the elements of $\Rchk$-matrices, so we have to
study the expansion of \eex{2.12} and\rx{2.13}. In accordance with
Proposition\rw{p2.1} we fix a number $M>0$ and study the action of
$\Rchk$-matrix only on those eigen-spaces $\vte\subset \vaivai$ for
which $\etzmnmax$.
 Therefore we assume that
%$\mmnzmn$
$0\leq m_1, m_2, n \leq (N-1)\amax$
in \eex{2.12}
and\rx{2.13}. As a result, the expansion of matrix elements is
achieved by substituting $1+\xh$ instead of $\spq$:
\qq
\Rchk(\fmot) & = &
\qfao \qmat \snz \,\plmn{\poh}
\label{3.1}\\
&&\qquad\times
\plnmc{\poh}\;
\spq^{-\a m_2}\, \poh^{m_2}\, \fmnoti
\nonumber
\qqq
and expanding the resulting expression in powers of $\xh$.

We left $\qa$ intact in the \rhs of \ex{3.1}. This gives us two
options for expansion. An evaluation of the \rhs in \ex{2.26}
requires us to expand in powers of $\xh$ while keeping $\a$ fixed.
This means that $1-\qma$ is of order $\xh$. Then the coefficients in
front of powers of $\xh$ in the expansion of the matrix
elements\rx{3.1} can be expanded themselves in $\xh$. This is not
dangerous for our purposes, because $1-\qma$ never appears in
denominator, so any further expansion goes in positive powers of
$\xh$ and does not change our estimates of `smallness'. The second
option is to keep $1-\qma$ constant and small as $\htz$ and use
\ex{3.1} to get at least a formal expansion in powers of $\xh$ and
$1-\qma$. We will excercise this option in the end of Section\rw{s5}
in order to prove the relation\rx{1.16}. Note that from the technical
point of view both options lead to the same expansion formulas. Only
the interpretation is different.

From now on we assume that $\a$ is kept constant as $\htz$. Let us
summarize the results of expanding the matrix elements\rx{3.1} in
powers of $\xh$. For every elementary positive braid in $B_N$ we
associate a parametrized $\Rchk$-matrix:
\qq
\Rpchke(\fmot) & = &
\snz \plmonf \lrbc{ e^{\e_{12}}a }^{n} \lrbc{e^{\e_1} }^{m_1}
\lrbc{ e^{e_2} \qma }^{m_2}
\nonumber\\
&&\qquad\qquad\qquad\qquad\qquad\qquad
\times\;
\fmnoti.
\label{3.2}
\qqq
This $\Rpchk$-matrix depends on parameters $a,\e_1,\e_2,\e_{12}$
`attached' to the elementary braid. For a negative elementary braid
we associate the `inverse' matrix
\qq
\Rpchkie(\fmot) & = &
\snz \plmtnf
\lrbc{e^{\e_{12}}a\p }^n \lrbc{ e^{\e_1} \qa }^{m_1}
\lrbc{ e^{\e_2} }^{m_2}
\nonumber\\
&&\qquad\qquad\qquad\qquad\qquad\qquad
\times\;
\fmnotii.
\label{3.3}
\qqq
Let us denote by $\spar$ the sets of parameters $a$, $a\p$ and
$\e_1,\e_2,\e_{12}$ associated to all elementary braids in a
particular expression of $B_N$. We construct a parametrized braid
operator $\Bhnpe$ acting on $\vaitn$, from the matrices $\Rpchk$ and
$\Rpchki$ in the same way as we constructed $\Bhn$ from $\Rchk$ and
$\Rchk^{-1}$. Since \eex{3.2} and\rx{3.3} do not constitute a
homomorphism of the braid group, the operator $\Bhnp$ depends on the
choice of a presentation of $B_N$ as a product of elementary braids.
%%%%%%%%%%%%%%%%%
\begin{proposition}
\label{p3.1}
For a given braid $B_N$
and its presentation as a product of elementary braids
there exists a set of polynomials
$$\Tje,\qquad \jgo$$
such that
\qq
\deg T_j \leq 2j
\label{3.4}
\qqq
and for a fixed $\eta$ and for any positive $M$
\qq
\lefteqn{
\qateb\, \qatn\, \qmet\, \Trfovnoe \Bhn
}
\label{3.5}\\
&=& \qhane \ffd \fsd
\nonumber\\
&&\qquad\times
\left.
\keta \Trfovnoe \Bhnpe \bckae + \ohmo.
\nonumber
\qqq
In our notations here $\{\cdot \}$ means all elements of the set.
\end{proposition}
%%%%%%%%%%%%%%%%%
\pr{Proposition}{p3.1}
Let us assume that the following lemma is true (see Appendix for the
proof):
%%%%%%%%%%%%%%%%%%%%%%%%
\begin{lemma}
\label{l3.1}
There exist two sets of polynomials
\qq
\Tjkro(m,n),\Tjrt(m,n)\in\Qmn, \qquad
\deg \Tjkro \leq j+k, \qquad \deg \Tjrt \leq 2j
\label{3.6}
\qqq
such that
the expansion of \ex{3.1} in powers of $\xh$ can be presented as
\qq
\Rchk(\fmot) & = &
\qfao \, \qmat \!\snz \!\!\left( \fone{n}{m_2}{m_2+n} \right)
\nonumber\\
&&\times\,
\lrbc{\ftwo{m_1}{n}} \lrbc{\fthree{m_2}}
\nonumber\\
&&
\qquad
\times\,
\plmonf\, (1-\qma)^n\,\spq^{-\a m_2} \, \fmnoti.
\label{3.7}
\qqq
\end{lemma}
A formula for $\Rchk^{-1}$ can be obtained from \ex{3.7} by
conjugation with $P$ (\ie, by the substitution\rx{2.1012}) and
substitutions $\spq\rightarrow \spq^{-1}$ and
$\xh\rightarrow \sujo (-1)^j\xh^j$.

%%%%%%%%%%%%%%%%%%%%%%%%%%%%%%%%%%%
For the practical purpose of using our formulas for actual
computation of the invariant polynomials $\vnkz$ of specific knots we
present the first polynomials $\Tjkro$ and $\Tjrt$:
\qq
\Trro_{1,0}(m_1,m_2)& = & - {1\over 2}(m_1+m_2+1)
\\
\Trro_{1,1}(m_1,m_2)& = & {1\over 6}
(m_1^2 + m_1 m_2 + m_2^2 + 3m_1 + 3m_2 + 2)
\\
\Trro_{2,0}(m_1,m_2)& = & {1\over 24}
(3m_1^2 + 6m_1m_2 + 3m_2^2 + 5m_2 + 7m_1 + 2)
\\
\Trt_1(m,n) & = & {1\over 2}(mn - n^2)
\\
\Trt_2(m,n) & = & {1\over 24}
(3 m^2 n^2 - 6 m n^3 + 3 n^4 + m^2 n - m n^2 - 5 m n + 5 n^2)
\qqq

The variables $m_1,m_2,n$ in the polynomials
\qq
\Tjkro(m_2,m_2+n),\;\Tjrt(m_1,n),\;\plzj(m_2-l)
\label{3.07}
\qqq
can be switched to derivatives
$\partial_{\e_1}$, $\partial_{\e_2}$, $\partial_{\e_{12}}$ with the
help of relations
\qq
\left.m_1 = \partial_{\e_1} e^{\e_1m_1}\right|_{\e_1=0}, \qquad
\left.m_2 = \partial_{\e_2} e^{\e_2m_2}\right|_{\e_2=0}, \qquad
\left. n = \partial_{\e_{12}} e^{\e_{12}n} \right|_{\e_{12}=0}.
\label{3.8}
\qqq
Another useful formula is
\qq
\left.
{\plzj (n-l)\over (1-\qma)^j} (1-\qma)^n =
\partial_a^j a^n \right|_{a=1-\qma}.
\label{3.9}
\qqq
A combination of \eex{3.8} and\rx{3.9} allows us to rewrite \ex{3.7}.
%%%%%%%%%%%%%%%%%%%%%
\begin{corollary}
\label{c3.1}
The expansion\rx{3.7} can be expressed in terms of the parametrized
matrix $\Rpchk$:
\qq
\Rchk
& = &
\qfao \, \qmat \left(
 1+
\sujo\xh^j\,\partial_a^j
     \skz
\xh^k\,\Tjkro(\partial_{\e_2},\partial_{\e_2}+\partial_{\e_{12}} )
%
%\partial_{\e_{12}}}{\partial_{\e_2}}{\partial_{\e_2} +
%\partial_{\e_{12}}}
%
\right)
\nonumber\\
&&\times\,
\lrbc{\ftwo{\partial_{\e_1}}{\partial_{\e_{12}}}}
\lrbc{\fthree{\partial_{\e_2}}}
\nonumber\\
&&
\left.
\qquad\times\,
\Rpchke \bcae.
\label{3.10}
\qqq
\end{corollary}
%%%%%%%%%%%%%%%%%%%%%%%%
The formula\rx{3.5} follows easily from this corollary. The
polynomials $T_j$ are formed by the products of
polynomials\rx{3.07} coming from elementary positive braids and by
the products of their counterparts coming from elementary negative
braids. The factors $\qfao$ of individual $\Rchk$-matrices\rx{3.10}
are canceled by the factor $\qateb$ in the \lhs of \ex{3.5}. The
factors $\qmat$ of individual $\Rchk$-matrices combine together
with the factor $\qatn$ from the \lhs of \ex{3.5} into the single
factor $\qhane$ in the \rhs of that equation. The factor $\qmet$ in
the \lhs of \ex{3.5} is represented by the factor
\qq
\left.
\fsd \keta \right|_{\kappa=1}
\label{3.11}
\qqq
in its \rhs. \qed

%%%%%%%%%%%%%%%%%%%%%%%%%%%%%%%%%%%%%%%%%%%
\nsection{The trace of parametrized $\Rchk$-matrix}
\label{s4}

\hyphenation{pa-ra-met-riz-ed}
The Proposition\rw{p3.1} presents the braid operator $\Bhn$ as a sum
of derivatives of the parametrized braid operator $\Bhnp$. Therefore
calculating the sum
\qq
\setz \keta \Trfovnoe \Bhnp
\label{4.1}
\qqq
is a natural step towards evaluation of the \rhs of \ex{2.26}.

The sum\rx{4.1} can be expressed in terms of the action of $\Bhnp$ on
the subspace $\vno\subset\vaitn$. We denote this action as $\Bhtnp$.
We express $\Bhtnp$ as a product of matrices $\Rpchk$ and $\Rpchki$
acting on $\vno$.

The natural basis in the space $\vno$ is formed by the vectors
\qq
e_j = f_0^{\otimes(j-1)} \otimes f_1 \otimes f_0^{\otimes(N-j)},
\qquad\ojN.
\label{4.2}
\qqq
The matrix $\Rpchk$ corresponding to the elementary positive braid
$\sjo$ transforms only the vectors $e_j$ and $e_{j+1}$. This
transformation is presented by a $2\times 2$ matrix
\qq
\Rhpchke =
\lrbc{
\begin{array}{cc}
e^{\e_1+\e_{12}}a & e^{\e_2}\qma \\
e^{\e_1} & 0
\end{array}
}.
\label{4.3}
\qqq
The matrix $\Rhpchki$ corresponding to the negative braid $\sjo^{-1}$
acts on the vectors $e_j$, $e_{j+1}$ as
\qq
\Rhpchkie =
\lrbc{
\begin{array}{cc}
0 & e^{\e_2}\\
e^{\e_1} \qa & e^{\e_2 + \e_{12}} a\p
\end{array}
}.
\label{4.4}
\qqq
Thus the matrix $\Bhtnp$ is a product of $N\times N$ matrices
containing $2\times 2$ blocks\rx{4.3} and\rx{4.4}.

%%%%%%%%%%%%%%%%%%%%%%%%%%
\begin{proposition}
\label{p4.1}
For a braid $B_N$ whose closure is a knot, there exist the constants
\qq
\delta_\e,\delta_a,\delta_{a\p},\delta>0
\qqq
such that if
\qq
|\{\e\}| < \delta_\e, \qquad
|\{a\}- (1-\qma)| < \delta_a, \qquad
|\{a\p\} - (1-\qa)|<\delta_{a\p}, \qquad
|1-\qa| < \delta
\qqq
(in our notations here $\{\cdot\}$ means `for any element of the
set') and $|\lambda|\leq 1$, then
\qq
\setz \leta \Trfovnoe \Bhnp =
{1\over
\det_{\fovnomo} \lrbc{1 - \lambda \Bhtnp} },
\label{4.5}
\qqq
here $\det_{\fovnomo}$ means that the matrix $\Bhtnp$ is projected
onto a subspace $\fovnomo\subset\vno$ spanned by the vectors $e_j$,
$\tjN$, along the vector $e_1$, and the determinant is calculated in
that subspace.
\end{proposition}
%%%%%%%%%%%%%%%%%%%
\pr{Proposition}{p4.1}
The formula\rx{4.5} follows from the fact that the operator $\Bhtnp$
acting on $\vaitn$ is an endomorphism of a symmetric tensor algebra
generated by a linear transformation of the subspace $\vno$.

The choice of basis $f_m$, $m\geq 0$ for the spaces $\vai$ allows us
to identify the tensor product $\vaitn$ with the algebra of
polynomials $\Czn$. The identification is achieved by mapping the
basis vectors $\fmtn\in \vaitn$ to monomials $\zmn$. The subspaces
$\vne$ are mapped onto the polynomials of degree $\eta$.

Consider the action of parametrized $\Rchk$-matrix\rx{3.2} on the
algebra $\Czt$ which is isomorphic to $\vaivai$.
%%%%%%%%%%%%%%%%%%%%%%%%%
\begin{lemma}
\label{l4.1}
The action of $\Rpchk$-matrix\rx{3.2} on $\Czt$ is an endomorphism
which is generated by the linear transformation\rx{4.3} of the
variables $z_1$ and $z_2$. The action of $\Rpchki$ of \ex{3.2} is an
endomorphism generated by the matrix\rx{4.4} in the same way.
\end{lemma}
%%%%%%%%%%%%%%%%%%%%%%%%%%%%
\pr{Lemma}{l4.1}
We will check this lemma only for $\Rpchk$. Consider a basis monomial
$\zmot$ corresponding to the basis vector $\fmot$ of $\vaivai$. The
linear transformation\rx{4.3} turns it into a polynomial
\qq
\lefteqn{
\lrbc{ e^{\e_1+\e_{12}}a\,z_1 + e^{\e_1}\,z_2}^{m_1}
\lrbc{ e^{\e_2}\qma\,z_1}^{m_2}
}
\hspace*{2in}
\label{4.6}\\
& = &
\sznmo {m_1\choose n}
\lrbc{e^{\e_1+\e_{12}}a\,z_1}^{n}
\lrbc{e^{\e_1}\,z_2}^{m_1-n}
\lrbc{e^{\e_2}\qma\, z_1}^{m_2}
\nonumber\\
& = &
\snz \plmonf
\lrbc{e^{\e_{12}}a}^{n}
\lrbc{e^{\e_1}}^{m_1}
\lrbc{e^{\e_2}\qma}^{m_2}\,
\zmnoti.
\nonumber
\qqq
The coefficients of the polynomial in the \rhs of this equation match
the matrix elements of the $\Rpchk$-matrix\rx{3.2}. \qed

The parametrized matrix $\Bhnp$ is a product of matrices $\Rpchk$ and
$\Rpchki$, so it is also an endomorphism of $\Czn$ whose action is
determined by the linear transformation $\Bhtnp$ of the space of
variables $\zn$. A projection of $\Bhn$ onto a subspace
$\fovi\subset\vaitn$ is an endomorphism of subalgebra
$\Cztn\subset\Czn$. This endomorphism is determined by the projection
of $\Bhtnp$ acting on the space $\vno$ of $N$ variables $\zn$, onto
the subspace $\fovnomo$ of $(N-1)$ variables $\ztn$. Therefore the
following well-known lemma proves \ex{4.5}:
%%%%%%%%%%%%%%%%%%%%%%%%%%%%%%%%%%%
\begin{lemma}
\label{l4.2}
Let $\cO$ be an endomorphism of the algebra $\Czn$ generated by its
linear action $\tcO$ on the space $\vno$ of variables $\zn$. Suppose
that for some $\lambda\in \IC$ the absolute values of eigen-values of
the operator $\lambda \tcO$ are smaller than 1. Denote by $\vne$ the
space of polynomials of degree $\eta$. Then
\qq
\setz \leta \Trvne \cO = {1\over \detvno \lrbc{1-\lambda\tcO} }.
\label{4.7}
\qqq
\end{lemma}
%%%%%%%%%%%%%%%%%%%%%%%%

An easy way to prove \ex{4.7} is to note that if $\tcO$ is
diagonalizable then \ex{4.7} becomes the formula for the sum of a
convergent geometric series. Since diagonalizable matrices form a
dense set in the space of all matrices, this proves the lemma.

It remains to check the applicability of the Lemma\rw{l4.2} to the
formula\rx{4.5}. We have to verify that the eigenvalues of the
operator $\Bhtnp$ are small enough. Note that the operator
\qq
\left.
\Bhtnpe
\bcaae
\label{4.8}
\qqq
acts on the space $\vno$ by permuting the basis vectors $\eon$ by the
permutation corresponding to the braid $B_N$. Since the closure of
$B_N$ is a knot, this means that for any vector $e_j$, $\tjN$ some
power of the operator\rx{4.8} will map it to $e_1$. Therefore the
projection of the operator\rx{4.8} onto the space $\vno$ is
nilpotent. Thus all the eigenvalues of this projection are zero.
Therefore the eigenvalues of the operator $\Bhtnpe$ are small as long
as $\{e\}$ are small, $\{a\}$ are close to $1-\qma$, $\{a\p\}$ are
close to $1-\qa$ and $1-\qa$ is kept small. \qed

\nsection{The colored Jones polynomial and perturbed Burau
representation}
\label{s5}

The Propositions\rw{2.1},\rw{3.1} and\rw{4.1} bring us close to the
following lemma:
%%%%%%%%%%%%%%%%%%%%%
\begin{lemma}
\label{l5.1}
If a knot $\cK\subset S^3$ is presented as a closure of a braid
$B_N$, then for any $M>0$
\qq
\scomb \Dmnk\,\a^{2m}\,\xh^n & = &
\qhane \ffd
\nonumber\\
&&\qquad\times\,\fsd
\label{5.1}\\
&&\qquad\times\,
\left.
{1\over \detfovnomo \lrbc{ 1 - \kappa\Bhtnpe } }
\bckae \!\!\!\!\!\!+ \ohmo.
\nonumber
\qqq
\end{lemma}
%%%%%%%%%%%%%%%%%%%%%%%%%%%%%%%%%

\pr{Lemma}{l5.1}
The only problem in combining \eex{2.26},\rx{3.5} and\rx{4.5} is that
in \ex{2.26} the sum over $\eta$ goes up to
$(N-1)\amax$ while in \ex{4.5}
the sum goes up to infinity.
Therefore the lemma would follow from the following estimate: for any
$M>0$
\qq
&&\ffd\fsd
\label{y.1}\\
&&\qquad\times\,
\left.
\setgmnmax \keta \Trfovnoe \Bhtnpe \bckae = \ohmo.
\nonumber
\qqq
%
%One has to substitute $\amax$ instead of $M$ in this estimate in
%order to apply it to the combination of \eex{2.26},\rx{3.5}
%and\rx{4.5}.

The proof of \ex{y.1} requires the following lemma:
%%%%%%%%%%%%%%%%%%%%%%%%%%%%%%%%%%%
\begin{lemma}
\label{l5.2}
For any $M>0$ the sum
\qq
\setgmn \keta \Trfovnoe \Bhtnpe
\label{5.2}
\qqq
is of combined order $M+1$ in variables $\{a\}$ and $\{a\p\}$.
\end{lemma}
%%%%%%%%%%%%%%%%%%%%%%

\pr{Lemma}{l5.2}
In view of \ex{4.5},
\qq
\setgmn \keta \Trfovnoe \Bhtnp
&=& {1\over \detfovnomo \lrbc{ 1 - \kappa\Bhtnp } }
\label{A2.1}\\
&&\qquad\left.
-
\setzmn {1\over n!}\,\partial_\lambda^n\;
{1\over \detfovnomo \lrbc{ 1 - \kappa\lambda\Bhtnp }}
\bclz.
\nonumber
\qqq
The \rhs of this equation is the error term in Taylor series
expansion of the value of
\qq
{1\over \detfovnomo \lrbc{ 1 - \kappa\lambda\Bhtnp }}
\label{A2.01}
\qqq
at $\lambda=1$. Therefore there exists $\lambda_0\in [0,1]$ such that
\qq
\setgmn \keta \Trfovnoe \Bhtnp =
\left.
{1\over (\nom+1)! } \, \partial_\lambda^{\nom+1}\;
{1\over \detfovnomo \lrbc{ 1 - \kappa\lambda\Bhtnp }}\bcllz
\label{A2.2}
\qqq

The determinant in the denominator of \rx{A2.01} is a polynomial in
$\kappa\lambda$ of degree
$\dim \fovnomo = N-1$:
\qq
\detfovnomo \lrbc{ 1 - \kappa\lambda\Bhtnpe } =
1 + \sojNo (\kappa\lambda)^j\,\tQje,
\label{A2.3}
\qqq
here $\tQje$ are some functions. The matrix $\Bhtnp$ is formed by the
products of matrices $\Rhpchk$ and $\Rhpchki$. The matrix elements
of\rx{4.3} and\rx{4.4} are products of variables
$a$, $a\p$, $e^{\e_1}$, $e^{\e_2}$, $e^{\e_{12}}$, $\qa$, $\qma$. As
a result the matrix elements of $\Bhtnp$ are integer polynomials of
these variables and
\qq
\tQj \in \zzlong.
\label{A2.03}
\qqq

There is an anti-symmetric counterpart of the Lemma\rw{l4.2}. Let
$\cO$ be an operator acting in the $N$-dimensional space $V_N$. Then
\qq
\det(1-\lambda\cO) = \szjN (-1)^j\lambda^j \Tr_{\bigwedge^j V_N} \cO,
\label{A2.4}
\qqq
here $\bigwedge^j V_N$ is the anti-symmetric part of $V_N^{\otimes j}$
and the action of $\cO$ is canonically defined there.
It follows from \ex{A2.4} that
\qq
\tQje = (-1)^j \Tr_{\bigwedge^j(\fovnomo)} \Bhtnpe.
\label{A2.5}
\qqq
We claim that the \lhs of this equation is of combined order 1 in
$\{a\}$ and $\{a\p\}$. Indeed, \eex{4.3} and\rx{4.4} indicate that
the matrix elements of the action of the operator
\qq
\left.
\Bhtnpe
\bcaaz
\label{A2.6}
\qqq
on the basis vectors $e_j$ of \ex{4.2} are proportional to the matrix
elements of the permutation operator corresponding to the braid
$B_N$. Therefore diagonal elements of the operator\rx{A2.6}
with respect to the basis vectors
$e_{j_1}\wedge\cdots\wedge e_{j_m}$, $m\leq N-1$
should be equal to zero (otherwise
a set $\{j_1,\ldots,j_m\}$, $m\leq N-1$ would be invariant under the
permutation of $B_N$ and the closure of $B_N$ would be a link rather
than a knot). Thus the value of the \rhs of \ex{A2.5} at
$\{a\}= \{a\p\}=0$ is zero. Hence the polynomials $\tQj$ are of
combined order at least 1 in $\{a\}$ and $\{a\p\}$.

An application of $\nom+1$ derivatives over $\lambda$ to the fraction
\qq
{1\over \detfovnomo \lrbc{ 1 - \kappa\lambda\Bhtnpe }} =
{1\over 1 + \sojNo (\kappa\lambda)^j\,\tQje }
\label{A2.7}
\qqq
brings the product of at least $N$ polynomials $\tQj$ to the
numerator. Since each of them is of combined order 1 in $\{a\}$ and
$\{a\p\}$, this proves the Lemma\rw{l5.2}.\qed

We substitute $\amax$ instead of $M$ in the expression\rx{5.2} in
order to apply Lemma\rw{l5.2} to the proof of the estimate\rx{y.1}.
The expression in the \lhs of \ex{y.1} contains derivatives over
$\{a\}$ and $\{a\p\}$ which reduce the powers of $a$ and $a\p$
in\rx{5.2}. However, according to \ex{3.10}, each derivative over $a$
and $a\p$ is accompanied by a power of $\xh$ in front of it.
Ultimately we set
$\{a\}=1-\qma$, $\{a\p\} = 1-\qa$ which makes them small of order
$\xh$. Since $\amax \geq M$, this proves the estimate\rx{y.1} and
completes the proof of \ex{5.1}. \qed

%Thus we see that the extra `would be' terms for $\etgmn$ in
%the \rhs of \ex{2.26} are of order $\xh^{M+1}$ after the
%substitution\rx{3.5}. This proves \ex{5.1}. \qed

Now we are ready to prove the following Proposition:
%%%%%%%%%%%%%%%%%%%%%%%%%%%%%%%%%%%%%%
\begin{proposition}
\label{p5.1}
For a knot $\cK\subset S^3$ there exist polynomials
$\pnqqa\in \Qqa$ such that
\qq
\scomb \Dmnk\,\a^{2m}\,\xh^n =
\sznM\xh^n {P_n( (1+h)^\a, (1+h)^{-\a})\over \atnka{\ohat} }
+ \ohmo.
\label{5.3}
\qqq
\end{proposition}
%%%%%%%%%%%%%%%%%%%%%%%%%%%

\pr{Proposition}{p5.1}
The key observation is that the matrices
\qq
\left.
\Rhpchke\bcaez & = &
\lrbc{
\begin{array}{cc}
1-\qma & \qma \\ 1 & 0
\end{array} },
\label{5.4}\\
\left.
\Rhpchkie\bcapez & = &
\lrbc{
\begin{array}{cc}
0 & 1 \\ \qa & 1-\qa
\end{array} }
\label{5.5}
\qqq
coincide with the Burau matrix and its inverse respectively if we set
\qq
t = \qma
\label{5.05}
\qqq
in the notations of\cx{KS}. Therefore the matrix
\qq
\left.
\Bhtnpe
\bcaape
\label{5.6}
\qqq
reproduces the Burau representation of the braid group and
\qq
\left.
\qahebn \det\nolimits_{\fovnomo}\!\! \lrbc{ 1 - \kappa \Bhtnpe }\bckae
= \aka{\qda}
\label{5.7}
\qqq
(see eq.~(5.22) of\cx{KS}).

The structure of the determinant in the \rhs of \ex{5.1} is described
by \eex{A2.3},\rx{A2.03} (we have to set $\lambda=1$ there). The
derivatives of the polynomial $\Tjed$, as well as the derivatives
$\partial_\kappa^j$, increase the power of the determinant in the
denominator and bring the polynomials $\tQje$ and their derivatives
to the numerator. The bound\rx{3.4} limits the `extra' powers of the
determinant in the denominator by twice the power of $\xh$. Thus we
recover the structure of the \rhs of \ex{5.3}. The expansion of the
prefactor
\qq
\qahneb \ohhne
\qqq
in the \rhs of \ex{5.1} in powers of $\xh$ does not change this
structure because
\qq
N-1-e(B_N) \in 2\ZZ
\label{5.8}
\qqq
if the closure of $B_N$ is a knot. \qed

The Proposition\rw{p1.1} follows from the Proposition\rw{p5.1} and
the next two lemmas which describe the properties of the polynomials
$\pnqqa$.
%%%%%%%%%%%%%%%%%%%
\begin{lemma}
\label{l5.3}
The polynomials $\pnqqa$ can be reexpressed as
\qq
\pnk{\qda} \in \Qdas.
\label{5.9}
\qqq
\end{lemma}
%%%%%%%%%%%%%%%%%%%%%%%
\begin{lemma}
\label{l5.4}
The polynomials $\pnqqa$ have integer coefficients:
\qq
\pnqqa \in \zqqa.
\label{5.10}
\qqq
\end{lemma}
%%%%%%%%%%%%%%%%%%%%%%%%%%%

\pr{Lemma}{l5.3}
The \lhs of \ex{5.3} is an even function of $\a$. The
Alexander-Conway polynomial of a knot
$\aka{\qda}$ is also an even function of $\a$. Therefore the
polynomials $\pnqqa$ are even in $\a$. Since for any $n\in\ZZ$,
\qq
\spq^{n\a} + \spq^{-n\a} = \Zdas,
\label{5.11}
\qqq
this proves the Lemma\rw{l5.3}.\qed

\pr{Lemma}{l5.4}
Let us use \ex{3.1} in order to expand the matrix elements of $\Rchk$
in powers of $\xh$ as we keep $1-\qma$ small but fixed. We will apply
the resulting formula towards the similar expansion of the trace
\qq
\tebn = \qateb\,\qatn \keta \detfovnoe \Bhn
\label{5.12}
\qqq
for the fixed value of $\eta$. In other words, we want to get an
expansion
\qq
\tebn = \snz \xh^n\,\tene.
\label{5.13}
\qqq

Let us first check that
\qq
\tene \in \zqqa.
\label{5.14}
\qqq
Indeed, as we know, the factors $\qfao\,\qmat$ of the
matrices\rx{5.7} combine with the factor $\qateb\,\qatn$ into a
single factor $\qhane$. Then, according to\rx{5.8},
\qq
\qhane = \qahneb\,\ohhne \in \zqqah.
\label{5.014}
\qqq
As for the other parts of the matrix elements\rx{3.1}, the relation
\qq
\plnmc{(1+\xh)} \in \zh
\label{5.15}
\qqq
shows that
\qq
\lefteqn{
\plmn{(1+\xh)}\,\plnmc{(1+\xh)}
}\hspace*{3in}
\nonumber\\
&&\times\,\spq^{-\a m_2} (1+\xh)^{m_2}
\in \zqqah.
\label{5.16}
\qqq
A similar relation holds for the elements of the inverse matrix
$\Rchk^{-1}$. The remaining piece of the \rhs of \ex{5.12}
\qq
\keta = (1+\xh)^{-\eta} \in \zhh
\label{5.17}
\qqq
presents no problem. Thus we proved relation\rx{5.14}.

The Propositions\rw{p3.1},\rw{p4.1} allow us to construct a
generating function of the traces\rx{5.12} as a formal power series
in $\xh$. More precisely, if we define a set of functions $\vnee$ by
the relation
\qq
\lefteqn{
\snz \xh^n\, \vnee
}
\label{5.19}
\\
& = &
\qahneb\,\ohhne \ffd
\nonumber\\
&&\qquad\times\,
\fsd
\nonumber\\
&&\qquad\times\,
\left.
{1\over \detfovnomo \lrbc{1 - \kappa\lambda \Bhtnpe} }
\bckae.
\qqq
then they generate the traces $\ten$ by the formula
\qq
\vnee = \setz \leta\, \tene,
\label{5.18}
\qqq

The structure of the determinant in the denominator of the \rhs of
\ex{5.19} is described again by \eex{A2.3},\rx{A2.03}. Therefore,
similar to the proof of the Proposition\rw{p5.1}, we observe that the
derivatives in the \lhs of \ex{5.18} add to the power of denominator
and bring the polynomials $\tQj$ and their derivatives to the
numerator. Thus we conclude that there exist the polynomials
\qq
\qnke\in\Qqa
\label{5.20}
\qqq
such that
\qq
\vnee = { \szktn \lambda^k\,\qnke \over
\lrbc{ 1 + \sojNo \lambda^j\,\qje} ^{2n+1} },
\label{5.21}
\qqq
here
\qq
\left.
\qje = \tQje\bcaape \in \zqqa.
\label{5.021}
\qqq
The bound $\zktn$ on the power of $\lambda$ in the numerator comes
from the fact that the maximum power of $\lambda$ in denominator is
$N-1$, each derivative brings one term $\lambda^j\,\tQj$ to the
numerator and the maximum number of derivatives at $\xh^n$ in the
\rhs of \ex{5.19} is $2n$.

The generating function $\vnee$ of \ex{5.18} satisfies the property
that
\qq
\left.
\fpar\,\vnee\bclz = \tene.
\label{5.22}
\qqq
This equation together with the relation\rx{5.14} and \ex{5.21}
implies that for all $\eta\geq 0$
\qq
\left.
\fpar\,
{ \szktn \lambda^k\,\qnke \over
\lrbc{ 1 + \sojNo \lambda^j\,\qje} ^{2n+1} }
\bclz \in \zqqa.
\label{5.23}
\qqq

A combination of relations\rx{5.021} and\rx{5.23} means that
\qq
\qnke \in \zqqa.
\label{5.24}
\qqq
Indeed, for a given $n$ let $k_0$ be the smallest value of $k$ for
which\rx{5.24} is not true. Then
%it would follow from \eex{5.21}
%and\rx{5.22} that
the relation\rx{5.23} is not satisfied for
%$T_{k_0,n}$.
$\eta = k_0$.

Comparing \eex{5.1},\rx{5.3},\rx{5.19} and\rx{5.21}, we see that
\qq
\pnqqa = \szktn \qnke.
\label{5.25}
\qqq
This together with the relation\rx{5.24} proves the
Lemma\rw{l5.4}.\qed

This completes the proof of Proposition\rw{p1.1}.

\nsection{Discussion}

The method of calculation of the coefficients $\vnkz$ in the
expansion\rx{1.13} of the colored Jones polynomial, which we
suggested in this paper, requires a presentation of the knot $\cK$ as
a closure of a braid. However, there are $R$-matrix formulas for the
colored Jones polynomial of a general blackboard picutre of a knot. A
direct calculation of the coefficients in this case requires a
certain generalization. The present method is based on a formula for
the trace of an endomorphism of a symmetric algebra. In case of a
blackboard picture, the number of strands is changing as one moves
upwards. Our calculation can be adapted for this case if we use the
so-called `holomorphic representation' of the homomorphisms of
symmetric algebras. This representation presents the homomorphisms as
integral operators with gaussian kernels acting on the algebra of
polynomials.

The Melvin-Morton expansion of the colored Jones polynomial comes in
the holomorphic representation approach from the stationary phase
approximation applied to a certain finite-dimensional integral. We
find this setting a bit more universal than taking derivatives of an
inverse determinant, as we have done here.

A slight modification of the procedure which led to the \rhs of
\ex{5.1}, can be applied to links. From the purely
technical point of view, this is an almost straightforward
exercise. This calculation produces invariants of links. We expect
these invariants to reproduce the Melvin-Morton expansion of the
colored Jones polynomial of a link only after a special
`step-by-step' procedure.

%However, we do not expect these new invariants to match
%the Melvin-Morton expansion of the Jones polynomial of a link
%directly.

The relation between the new invariants of links and the colored
Jones polynomial can be explained (at the level of a conjecture) with
the help of path integral arguments.
%We need to introduce some path integral arguments in order to explain
%our expectations in case of links.
Let
\qq
\spq = e^{2\pi i\over K}, \qquad K\in\ZZ.
\qqq
Consider the colored polynomial of a link in the limit of
$K\rightarrow \infty$, while we keep the ratios
\qq
a_j = {\a_j\over K}
\label{6.01}
\qqq
($\a_j$ being the colors of the link components) constant. According
to\cx{EMSS}, the Jones polynomial in this limit can be presented as a
path integral of the exponential of the Chern-Simons action. The
integration should go over all $SU(2)$ connections $A_\mu$ in the
link complement which satisfy the boundary condition: up to a
conjugation, in the fundamental representation of $SU(2)$
\qq
\PexpA{\cC_j} = \exp \lrbs{
i\pi a_j \lrbc{
\begin{array}{cc}1&0\\0&-1 \end{array} } },
\label{6.1}
\qqq
here $\PexpA{\cC_j}$ is a physical notation for the holonomies of the
connection $A_\mu$ along the meridians $\cC_j$ of the link
components.

In the limit of $K\rightarrow\infty$ the path integral can be
calculated in the stationary phase approximation. The stationary
points are flat connections in the link complement which satisfy the
boundary conditions\rx{6.1}. In case of a knot, if $a={\a\over K}$ is
small enough, then there is only one flat connection which
satisfies\rx{6.1}. Its contribution generates the expansion\rx{1.13}.
This connection is reducible, all its holonomies lie in the same
subgroup $U(1)\subset SU(2)$. This reducibility corresponds to the
fact that the expansion\rx{1.13} came in our calculations from the
contribution of the `edge' of the space $\vatn$. Within the path
integral interpretation, the `top' of $\vatn$ corresponds to the
situation when all the holonimies are aligned along the same subgroup
$U(1)\subset SU(2)$.

In case of a link, generally there are flat irreducible $SU(2)$
connections in the link complement satisfying the condition\rx{6.1}
even if the phases $a_j$ are arbitrarily small. Therefore we do not
expect that the application of \ex{5.1} to a braid whose closure is a
link, would yield the actual expansion of the colored Jones
polynomial of that link. We conjecture that the \rhs of \ex{5.1}
gives an expansion of the contribution of a reducible connection to
the Jones polynomial. This contribution was discussed in\cx{Ro3} in
the context of the Reshetikhin formula for the Jones polynomial. It
seems that this contribution is an invariant of the link in itself.
We also think that it determines the perturbative invariants as
well as $p$-adic properties of the Witten-Reshetikhin-Turaev
invariant of rational homology spheres constructed by a surgery on
the link.

We will address the issues of the blackboard calculations,
holomorphic representation and new invariants of links in the future
publication.

\section*{Acknowledgements}

I am thankful to D.~Bar-Natan, J.~Birman, P.~Deligne, V.~Drinfeld,
D.~Freed, S.~Garoufalidis, L.~Kauffman, M.~Kontsevich, P.~Melvin,
R.~Nepomechie, D.~Thurston, V.~Turaev, A.~Vaintrob,
A.~Varchenko, Z.~Wang, E.Witten and C.~Zachos for many useful
discussions.

This work was supported by the National Science Foundation
under Grants No. PHY-92 09978 and DMS 9304580. Part of this work was
done during my visit to Argonne National Laboratory. I want to thank
Professor C.~Zachos for his hospitality.

%\nappendix{1}

\nappendixe

\pr{Lemma}{l3.1}
This lemma follows from the next two lemmas.
%%%%%%%%%%%%%%%%%%%%%%%%%%%
\begin{lemma}
\label{lA.1}
There exist the polynomials
$\Tjrt(m,n)$, $\deg\Tjrt \leq 2j$
such that
\qq
\plnmwc = \plnmnf\,
\lrbc{\ftwo{m}{n} }.
\label{A.1}
\qqq
\end{lemma}
%%%%%%%%%%%%%%%%%%%%%%%%

\pr{Lemma}{lA.1}
It is enough to prove that there exist the polynomials
$\tTjrt(n)$, $\deg \tTjrt \leq 2j$ such that
\qq
\plon {\ohlo\over \xh} = n!\,\lrbc{ 1 + \sujo \xh^j \tTjrt(n)}.
\label{A.2}
\qqq
Lemma\rw{lA.1} would follow because the \lhs of \ex{A.1} can be
presented as a combination of three products of the type\rx{A.2}.

Consider the logarithm of the \lhs of \ex{A.2}:
\qq
\log \lrbc{ \plon {\ohlo\over \xh} }
& = &
\log(n!) + \slon \log\lrbc{ 1+\sko \xh^k\,{\pmzk (l-m)\over (k+1)!} }
\label{A.3}\\
& = &
\log(n!) - \slon\skpo {(-1)^{k\p} \over k\p}
\lrbc{ \sko \xh^k {\pmzk (l-m)\over (k+1)!} }^{k\p}.
\nonumber
\qqq
The coefficient at a given power $\xh^j$ in the \rhs of \ex{A.3} is a
polynomial in $l$ of degree $j$. The sum $\slon$ turns it into a
polynomial in $n$ of degree $j+1$. Since $j\geq 1$, this degree is
not greater than $2j$. Taking the exponential of \ex{A.3} and
expanding it in powers of $\xh$ gives us the expansion\rx{A.1}. \qed

%%%%%%%%%%%%%%%%%%%%%%%%%%%%%%%%%
\begin{lemma}
\label{lA.2}
There exist the polynomials
$\Tjkro(m_1,m_2)$, $\deg\Tjkro\leq j+k$, such that
\qq
\plmpon \lrbc{ (1+\xh)^{-l} - \qma} & = &
\lrbc{ \fone{n}{m}{m+n} }
\nonumber\\
&&\qquad\times\,
(1-\qma)^n.
\label{A.4}
\qqq
\end{lemma}
%%%%%%%%%%%%%%%%%%%%%%%%%%%%%%%%%

\pr{Lemma}{lA.2}
We expand the \lhs of \ex{A.4} in powers of $\ohlmo$:
\qq
\plmpon \lrbc{ (1+\xh)^{-l} - \qma} & = &
\plmpon \lrbs{ \lrbc{\ohlmo} - \lrbc{1-\qma} }
\label{A.5}\\
& = &
(1-\qma)^n \sjz {\xh^j\over (1-\qma)^j }\, S_j(m, m+n;\xh).
\nonumber\\
S_j(m_1,m_2;\xh) & = & \sum_{\mollmt} \poij {(1+\xh)^{-l_i} - 1\over
\xh}.
\label{A.6}
\qqq
Now it remains to prove the following:
%%%%%%%%%%%%%%%%%%%%%%%%%%%%%%%%%
\begin{lemma}
\label{lA.3}
There exist the polynomials
$\Tjkro(m_1,m_2)$, $\deg\Tjkro\leq j+k$, such that
\qq
S_j(m_1,m_2;\xh) = \lrbc{ \pzijo (m_2-m_1-i) }
\skz \xh^k\,\Tjkro(m_1,m_2).
\label{A.7}
\qqq
\end{lemma}
%%%%%%%%%%%%%%%%%%%%%%%%%%%%%%%%

\pr{Lemma}{lA.3}
We prove this lemma by induction. When $j=0$, the lemma is obvious:
\qq
S_0(m_1,m_2;\xh) = 1,\qquad
T_{0,k}(m_1,m_2) = \delta_{0,k}.
\label{A.8}
\qqq

Suppose that the lemma is true for a particular value of $j$. The sum
$S_{j+1}(m_1,m_2;\xh)$ can be expressed in terms of $S_j$:
\qq
S_{j+1}(m_1,m_2;\xh) =
\smotl S_j(m_1,l-1;\xh)\, {\ohlmo \over \xh}.
\label{A.9}
\qqq
Expanding the last factor
\qq
{\ohlmo\over \xh} = \skz \xh^k\, {\pzik (l+i)\over (k+1)! }
\label{A.10}
\qqq
and substituting the formula\rx{A.7} for $S_j(m_1,l-1;\xh)$ we find
that
\qq
&
S_{j+1}(m_1,m_2;\xh)  =
\skz \xh^k\szkpk \Tojk(m_1,m_2),
\label{A.11}\\
&
\Tojk(m_1,m_2) = \smotl \poij (l-m_1-i)
\,{\pzikkp(l+i) \over (k-k\p+1)!}\,
\Tjkpro(m_1,l-1).
\label{A.12}
\qqq
The \rhs of \ex{A.12} contains finite sums of polynomials of
$m_1,l$, so it is clear that $\Tojk$ is a polynomial in $m_1,m_2$
such that
\qq
\deg\Tojk \leq 2(j+1) + k.
\label{A.13}
\qqq

The factor $\poij(l-m_1-i)$ in \ex{A.12} guarantees that
\qq
\Tojk(m_1,m_2) = 0 \qquad
\mbox{if $m_2-m_1\in \ZZ$, $1\leq m_2-m_1\leq j$.}
\label{A.14}
\qqq
Denote by $\Bnm$ a polynomial
\qq
\Bnm = \smmot m^n.
\label{A.15}
\qqq
Since $\Bnm=0$ for $m_1=m_2$, we conclude from \ex{A.12} that
\qq
\Tojk(m_1,m_2)=0 \qquad
\mbox{if $m_2-m_1=0$}.
\label{A.16}
\qqq
Equations\rx{A.14} and\rx{A.16} imply together that there exist the
polynomials $\tTojk(m_1,m_2)$ such that
\qq
\Tojk(m_1,m_2) = \lrbc{ \pzij (m_2-m_1 -l) } \,\tTojk(m_1,m_2).
\label{A.17}
\qqq
Then the polynomials
\qq
\Tjokro(m_1,m_2)=\szkpk \tTojk(m_1,m_2)
\label{A.18}
\qqq
satisfy \ex{A.7} at $j+1$. They also have the right degree in view of
\ex{A.13}. This proves Lemma\rw{lA.3}. \qed

\end{document}